\documentclass[12pt,a4paper]{article}

\usepackage{pdflscape}
\usepackage{amsmath}
\usepackage{amssymb}
\usepackage{latexsym}
\usepackage{srcltx}
\usepackage{graphics}
\usepackage{color}
\usepackage{epsfig}
\usepackage{color}

\usepackage{subcaption}

\usepackage{anyfontsize}
\usepackage{tabularx}

\usepackage{graphicx}
\usepackage{multirow}
\usepackage{array}

\usepackage{amssymb, amsmath, amsfonts, mathtools}
\usepackage{esint}
\usepackage{graphicx}
\usepackage{breqn}
\usepackage{float}
\usepackage{caption}
\usepackage{subcaption}
\usepackage{hhline}
\usepackage{multirow}
\usepackage{rotating}
\usepackage{fixltx2e}

\newcommand{\re}{{\mathbb R}}
\newcommand{\n}{{\mathbb N}}

\newcommand{\cA}{{\mathcal{A}}}

\newcommand{\cF}{{\mathcal{F}}}

\newcommand{\cB}{{\mathcal{B}}}

\newcommand{\cS}{{\mathcal{S}}}

\newcommand{\cM}{{\mathcal{M}}}

\newcommand{\bx}{{\boldsymbol x}}

\newcommand{\ba}{{\boldsymbol a}}
\newcommand{\bb}{{\boldsymbol b}}

\newcommand{\be}{{\boldsymbol e}}

\newcommand{\br}{{\boldsymbol r}}
\newcommand{\bv}{{\boldsymbol v}}
\newcommand{\bu}{{\boldsymbol u}}

\newcommand{\bh}{{\boldsymbol h}}

\newtheorem{theorem}{Theorem}
\newtheorem{prop}{Proposition}
\newtheorem{lemma}{Lemma}
\newtheorem{cor}{Corollary}
\newtheorem{remark}{Remark}
\newtheorem{ex}{Example}
\newtheorem{defi}{Definition}

\date{}

\author{Aleksandar Cvetkovi\' c, \thanks{GSSI (L'Aquila, Italy), }
Vladimir Yu.~Protasov
\thanks{University of L'Aquila (Italy), 
Moscow State University (Moscow, Russia),  {e-mail: \tt\small
v-protassov@yandex.ru}}}

\title{The greedy strategy for optimizing the Perron eigenvalue 
\thanks{
The research is supported by the RSF grant 20-11-20169
}}

\begin{document}
\maketitle

\begin{abstract}

We address the problems of  minimizing and of maximizing the spectral radius 
over a compact family of  non-negative matrices. Those problems being 
hard in general can be efficiently solved for some special  families. 
We consider the so-called product families, where each matrix is composed of 
rows chosen independently from given sets. A recently introduced  greedy method 
works very  fast. However, 
it is applicable mostly  for strictly positive matrices.  
For sparse matrices, it often diverges and gives a wrong answer.  
We present the ``selective greedy method'' that works equally well 
for all non-negative product families, including sparse ones. For this method, we prove 
a quadratic rate of convergence and  demonstrate its  efficiency in numerical examples.  
The numerical examples are realised for two cases: finite uncertainty sets and 
polyhedral uncertainty sets given by systems of linear inequalities. In dimensions up to 2000,  the 
matrices with minimal/maximal spectral radii in product families are 
found within a few iterations.   
Applications to  dynamical systems and to the graph theory are considered.

\bigskip

\smallskip

\begin{flushright}
\noindent  \textbf{AMS 2010} {\em subject
classification: 15B48, 90C26, 15A42}
\end{flushright}

\end{abstract}
\bigskip

\begin{center}
\large{\textbf{1. Introduction}}
\end{center}
\bigskip

The problem of minimizing or maximizing the spectral radius 
of a non-negative matrix over a general constrained set can be very hard. There are no efficient algorithms 
even for solving this problem  over a compact convex (say, polyhedral) set of 
positive matrices. 
This is because the objective function  is neither convex nor 
concave in matrix coefficients and, in general, non-Lipschitz. There 
may be  
many points of local extrema which are, moreover, hardly identified. Nevertheless, 
for some special sets of matrices efficient methods do exist. 
In this paper we consider  the so-called {\em product families}. 
Their  spectral properties were discovered and analysed 
by Blondel and Nesterov in~\cite{BN}; the first methods of optimizing the spectral radius 
over such families originated in~\cite{NP1}. One of them, the {\em spectral simplex method}, was further developed in~\cite{P16}. This method is quite simple in 
realization and has a fast convergence to the global minimum/maximum. 
Similar methods for optimising the spectral radius appeared in the literature in various contexts: in stochastic and entropy games~(see~\cite{Akian1} and references therein),  
 PageRanking~\cite{CJB, FABG}. The relation to the Markov Decision Processes 
 is discussed in detail in Remark~\ref{r.5} below.

A modification  of the spectral simplex method, 
the {\em greedy method}~\cite{Akian1, NP2}, demonstrates a 
very good efficiency.  Within 3-4 iterations,  
it finds the global minimum/maximum with a tight precision. 
Moreover, the number of iterations does not essentially grow with the dimension. However, for sparse matrices this method often suffers. It may either cycle or 
converge to a very rough estimate. We provide an  example of cycling in Section~3. 
In that example there is a product family of $3\times 3$ matrices
for which the 
 algorithm of greedy method cycles  and the value of the spectral radius stabilises
 far from the real optimal value. For large sparse matrices, this phenomenon occurs quite often. This is a serious disadvantage, since  most of 
applications deal with sparse matrices.  There is a standard perturbation technique, when all zeros 
of a sparse non-negative matrix are replaced  by small strictly positive entries. 
It is applied to avoid negative effects of sparsity: multiple Perron eigenvalues, 
divergence of the power method, etc. However, in our case it often 
 causes computational problems: 
 due to the high dimensions and the large number of matrices involved  we either obtain a rough approximation of the spectral radius or set the perturbation value too small, which leads to the cycling again.

In this paper we develop the greedy method, which works equally well for all 
product families of non-negative matrices, including sparse matrices. Numerical results demonstrated in Section~6 show that even in dimension of several thousands, the 
greedy method manages to find minimal and maximal spectral radii within a few iterations. We also estimate the rate of convergence and thus theoretically 
justify the efficiency of the method. Finally we consider several applications. 

Before introducing the main concepts let us define some notation. We 
denote the vectors by bold letters and their components by standard 
letters, so $\ba = (a_1, \ldots , a_d) \in \re^d$. 
The {\em support} of a non-negative vector is the set of indices of its 
nonzero components. For a non-negative 
$d\times d$ matrix~$A$, we denote by $a_{ij}$ its entries, by $\rho(A)$ its spectral radius, which is 
the maximal modulus of its eigenvalues. By the Perron-Frobenius theorem, 
 $\rho(A)$ is equal to the maximal non-negative eigenvalue $\lambda_{\max}$
 called the {\em leading eigenvalue} of the matrix~$A$. 
 The corresponding non-negative eigenvector~$\bv$ is also called {\em leading}.
 Note that the leading eigenvector may not be unique. 
 Each non-negative $n\times n$ matrix~$A$ is associated to 
 a directed graph with $n$ vertices numbered as $1, \ldots , n$ and defined as follows: there exists an edge from $i$ to $j$ if and only if $A_{ji} >0$. 
 A non-negative matrix is {\em reducible} if after some renumbering of coordinates it 
 gets a block upper-triangular form. 
The  matrix is irreducible precisely when 
 its graph is strongly connected. An irreducible matrix has a unique 
 up to normalization  leading eigenvector.    We usually normalize it as $\|\bv\| = 1$, where $\|\bv\| = \sqrt{(\bv, \bv)}$ is 
 the Euclidean norm. 
 
\begin{defi}\label{d.10}
Let $\cF_i \subset \re^d_+, \, i = 1, \ldots , d$, be arbitrary 
 nonempty compact sets 
referred to as {\em uncertainty sets}. 
Consider a matrix $A$ such that for each $i = 1, \ldots d$, 
the $i$th row of $A$ belongs to $\cF_i$.  The family of all those 
matrices~$A$ is denoted by $\cF$ and called a {\em product family} or a  
{\em family with product structure}. 
\end{defi}
Thus, we compose a matrix $A$ choosing each row from the corresponding uncertainty set. The family of all such matrices $\cF = \cF_1\times \cdots \times \cF_d$
is, in a sense, a product of uncertainty sets. Of course, not every compact set 
of non-negative matrices has a product structure. Nevertheless, the 
class of product families 
is important in many applications. For instance, in 
dynamical systems and asynchronous systems~\cite{Koz1, LA, P16},  
Markov Decision Processes~\cite{Al, Y}, PageRanking~\cite{FABG}, 
graph theory~\cite{ESS, Liu, NP1}, mathematical 
economics~\cite{CL, Log, BN}, game theory~\cite{Asarin, Akian1, HK, R}, 
matrix population models~(see~\cite{Log} and references therein), etc.   

We consider the problems of optimising the spectral radius
over a product family~$\cF$: 
\begin{equation}\label{eq.main}
\left\{
\begin{array}{l}
\rho(A)\ \to \ \min \, /\,  \max\\
A \, \in \, \cF
\end{array}
\right. 
\end{equation}
Both minimization and maximization 
problems can be efficiently solved~\cite{NP1}. The spectral simplex method
for finding global minimum and maximum demonstrates a fast convergence~\cite{P16}. For example, if all the  sets $\cF_i$ are two-element (so the set 
$\cF$ has a structure of a Boolean cube), then in dimension $d=100$, 
the matrices  with minimal and maximal spectral radii are  found 
(precisely!)
within $50 - 60$ iterations for 10-15 sec in a standard laptop. 
 Note that in this case $\cF$ contains $|\cF| = 2^{100} > 10^{30}$ matrices. 
If for the same dimension $d=100$, each uncertainty set consists of $100$
rows (so, $|\cF| = 10^{200}$), then the spectral simplex method 
performs about $200 - 300$ iteration and solves the problem~(\ref{eq.main})
in less than a minute. We write the algorithm in the next section, but its 
idea can be described within a few lines. Let us consider the maximization 
problem~(\ref{eq.main}). We start with an 
arbitrary matrix $A_1 \in \cF$ with rows $\ba_i \in \cF_i, \, i = 1, \ldots , d$, 
and compute its leading eigenvector $\bv_1$.
Then we solve the problem $(\bv_1, \ba)\to \max, \, \ba\in \cF_1$, in other words,  
we find an element from $\cF_1$ that makes the biggest projection onto the eigenvector $\bv_1$. If the row $\ba_1$ is already optimal, then we leave it and 
do the same with the second row and with the set~$\cF_2$, etc. 
If all rows of $A$ are optimal, then $A$ has the biggest spectral radius. 
Otherwise, we replace the row~$\ba_i$ by an optimal element from~$\cF_i$
and thus obtain the next matrix~$A_2$. Then we do the next iteration, etc. 
Thus, by one iteration we mean one computation of the 
leading eigenvector of a matrix, all other operations are cheap.

If all the sets~$\cF_i$ consist of strictly positive rows, then 
the spectral radius increases each iteration. 
So, we have a relaxation scheme which always converges to 
a global optimum. However, if the rows may have zeros, 
then the relaxation is non-strict (it may happen that $\rho(A_{k+1}) = \rho(A_k)$), 
and the algorithm may cycle or converge to a non-optimal matrix.  
In~\cite{P16} this trouble was resolved and the method was modified to be applicable 
 for sparse matrices. We will formulate the idea  in Theorem~A 
in the next section. 
\smallskip 

Recently in~\cite{Akian1, NP2} the spectral simplex method 
was further improved: in each iteration we maximize not one row but all 
rows simultaneously. According to numerical experiments, 
this slight modification leads to a significant improvement: 
in all practical examples the algorithm terminates within $3-4$ iterations.  
In the example above, with $d=100, \, |\cF_i|=100, i = 1, \ldots , d$, 
the absolute minimum/maximum is found within a few seconds and the number of iterations rarely exceeds $5$. Even for dimension $d=2000$, 
the number of iterations never exceeds~$10$. 
The authors of~\cite{Akian1} and of~\cite{NP2} came to problem~(\ref{eq.main})
from different angles. The paper~\cite{Akian1} studies problems from  entropy
games with a fixed number of states.  The paper~\cite{NP2} solves 
a different problem:  finding the closest stable matrix to a given matrix in the~$L_1$ norm. 
The authors derive a method similar to the policy iteration method and call it~{\em the greedy method}. 
We will use the latter name. This modification of the spectral simplex method, although has a much faster convergence, 
inherits the main disadvantage: 
it works  only for strictly positive matrices. For sparse matrices,  it often 
gets stuck, cycles, or converges to a wrong solution. 
The standard perturbation technique leads co computational problems
in dimensions, we discuss this issue in Section 3 and subsection 4.2.  
None of ideas from the work~\cite{P16} which successively modified  
the spectral simplex method helps  for the greedy method.  This is 
shown in example and discussed in  Section~3. 
In some of the litarature on entropy and stichastic games  (see~\cite{Akian1}) the positivity assumption  
is relaxed to irreducibility of all matrices~$A_k$ in each step, but this condition is sometimes  very restrictive and 
is never satisfied in practice for sparse matrices. In~\cite{NP2} 
the greedy method was modified for sparse matrices, but only 
for minimization problem and by a significant complication of the computational procedure. 
 Therefore, in this paper we attack the three main problems: 
\smallskip 

\textbf{The main problems}: 
\smallskip 

1. To extend the greedy method for sparse matrices while preserving its efficiency.  
\smallskip 

2. To estimate the rate of convergence of the greedy method. 
\smallskip 

3. To apply the greedy method to some of known problems. We consider two applications: 
optimisation of the spectral radius of graphs and finding the closest 
stable/unstable linear system.  
\smallskip 

\smallskip 

The first issue is solved in Section~4. We derive  the 
{\em selective greedy method} working for all kind of non-negative matrices. 
To this end we introduce and apply the notion of {\em selected leading eigenvector}. 
The new greedy method is as simple in realisation and as fast as the previous one. 

In Section~5 we prove two main theorems about the rate of convergence. 
We show that both the simplex and the greedy methods have a global linear 
rate of convergence. Then we show that the greedy method (but not the spectral simplex!) has a local quadratic convergence, which explains its efficiency. We also estimate the constants and parameters of the convergence. 

Numerical results for matrix families in various dimensions are reported in Section~6. 
Section~7 shows application of the greedy method to problems of dynamical systems and of the graph theory.  Finally, in Section~8, we 
discuss the details of practical implementation.

\begin{remark}\label{r.5}
{\em A relation of the spectral optimisation over product families to the 
 Markov Decision Processes (MDP) deserves a separate detailed consideration. 
The MDP provide a mathematical framework for decision making when outcomes are partly random and partly under the control of a decision maker. The theory of MDP is developed 
 since early 60s and is widely known and popular nowadays, see, for instance~\cite{Al, HK, P} and references therein. 
 
 Our spectral optimisation issue has a lot in common 
 with MDP, from the statement of the problem to methods of the solution. 
First, the MDP can  be described by means of  product families of non-negative matrices. Moreover, the 
{\em policy iteration algorithm}  for finding the optimal strategy in MDP looks very similar to the 
spectral simplex method and the greedy method. That algorithm  chooses the rows from uncertainty sets that are optimal 
in some sense (maximise the reward), then computes a special vector~$\bv$
(the value), then again chooses the rows using  this new vector, etc.  

Nevertheless, the problem of optimising the spectral radius over product family 
cannot be derived from MDP, and the  spectral simplex method and the greedy method 
are different from those existing in the MDP literature. The crucial difference is that MDP does not optimise 
the leading eigenvalue and  does not deal with the leading eigenvectors. Let us recall the 
main scheme of the policy iteration, the algorithm which seems to be   
similar to the spectral simplex and to the greedy method. For the sake of simplicity we describe  this algorithm in terms 
of product families. We are given a number~$\gamma \in [0,1]$ (discount factor) and finite uncertainty sets $\cF_i, \, i = 1, \ldots , d$. 
To each element $\ba_i^{(j)} \in \cF_i$, one associates a given non-negative row vector~$\br_i^{(j)}$ (the reward).   In each iteration of the algorithm, when 
we have a matrix~$\cA \in \cF$ and the corresponding rewards~$\br_1, \ldots , \br_d$, 
we find the vector~$\bv = (v_1, \ldots , v_d)\in \re^d$ by solving the following linear 
non-homogeneous system:  
\begin{equation}\label{eq.pia1}
v_i \ = \ (\ba_i, \bv)\ + \ \gamma\, (\ba_i, \br_i)\, , \quad i = 1, \ldots , d\, . 
\end{equation} 
 Then  we choose  new rows~$\ba_i' \in \cF_i, \, i = 1, \ldots , d$
 (with the corresponding~$\br_i'$) as a solution of the optimisation 
 problem: 
\begin{equation}\label{eq.pia2}
\max_{\ba_i' \in \cF_i}\ \Bigl[\, (\ba_i', \bv)\ + \ \gamma\, (\ba_i', \br_i') \, \Bigr]
\end{equation}  
 (now $\bv$ is fixed and the optimisation runs  over the sets~$\cF_i$). 
Then we find new~$\bv$ by solving linear system~(\ref{eq.pia1}) for the new matrix~$A'$, etc., 
until in some iteration we have~$A' = A$. In this case the matrix~$A$ is optimal, 
which corresponds to the optimal choice from the uncertainty sets. 
It is known that in case of finite uncertainty sets, the policy iteration algorithm converges 
within finite time. Moreover, that time is actually polynomial~\cite{Y}.  
Efficient estimates for the number or iterations in terms of $d, |\cF_i|,$ and $\gamma$
are obtained in~\cite{HMZ, MS, PY}. There are various versions of the policy iteration algorithm, see~\cite{GG, P}. 
The MDP problems can also be solved by different techniques, including linear 
programming~\cite{Y}. 

We see that the MDP problem is different from optimising the spectral radius, although 
it also deals with product families. 
The policy iteration, in spite of the similarity of the structure,  is different from both 
spectral simplex and greedy method. It solves a different problem and does not deal with leading 
 eigenvectors 
using another vector~$\bv$ obtained from linear equation~(\ref{eq.pia1}). 

A natural question arises if some of  MDP algorithms  and the greedy method are equivalent and can be derived from each other? 
To the best of our knowledge, the answer is negative. One of possible 
arguments is that for algorithms from MDP literature, in particular, for the policy iteration, there are guaranteed estimates for the rate of convergence, while  the greedy 
algorithm and the spectral simplex methods may not converge at all and may suffer of sparsity, cycling, 
etc., see~\cite{P16} and Section~3 of this paper. So, in our opinion, the results of this paper are hardly 
applicable to MDP and vice versa.  
}
\end{remark}

\medskip 

\textbf{Novelty}. Let us formulate the main novelty of our results.  
\smallskip 

 \textbf{1.} We develop the greedy method of optimising the 
spectral radius over a product family of non-negative matrices. This method works
 equally well for infinite uncertainty sets. For polyhedral uncertainty sets given by linear inequalities, where the number of vertices can be very large, the method finds 
 the precise optimal value within a few iterations (Section~6).  For all other 
 compact uncertainty sets, it finds approximate values.
\smallskip 

\textbf{2.} We prove the linear convergence for the spectral simplex method and the quadratic converges for the greedy method. The parameters of linear and quadratic convergence are estimates in terms of geometry of the uncertainty sets. 
\smallskip 

\textbf{3.} We derive a simple procedure to extend the greedy method to 
sparse matrices without loss of its efficiency. The extension is based on 
the use of {\em selective leading eigenvectors}.  
\smallskip

\bigskip 

\begin{center}
\large{\textbf{2. The minimizing/maximizing spectral radius \\
of product families}}
\end{center}
\medskip

\begin{center}
\textbf{2.1.  The description  of the methods}
\end{center}
\medskip

We have a product family~$\cA$ with compact uncertainty sets $\cF_1, \ldots , \cF_d$
which are supposed to be non-negative (consist of entrywise non-negative vectors).
The problem is to find the maximal  and the minimal spectral radius of matrices form~$\cA$.

A matrix $A \in \cA$ is said to be {\em minimal in each
row} if it has  a leading eigenvector $\bv \ge 0$ such that $(\ba_i, \bv) = \min_{\bb_i \in \cF_i}(\bb_i , \bv)$
for all $i = 1, \ldots , d$. It is {\em maximal in each row} if it has a leading eigenvector
$\bv > 0$ such that $(\ba_i, \bv) = \max_{\bb_i \in \cF_i}(\bb_i , \bv), \, i = 1, \ldots , d$.
Those notation are not completely analogous: the minimality in each row is defined with respect to an arbitrary  leading eigenvector~$\bv \ge 0$, while the maximality needs a strictly positive~$\bv$.
\smallskip 

We  first briefly describe the ideas of the methods and then 
write formal procedures. Optimisation of spectral radius of a product family is done by a relaxation scheme. We consider first the maximization problem. 

Starting with some matrix $A_1 \in \cF$
we build a sequences of matrices $A_1, A_2, \ldots $ such that 
$\rho(A_1) \le \rho(A_2) \le \cdots $ Every time the next matrix~$A_{k+1}$ 
is constructed from~$A_k$ by the same rule. The rule depends on the method:  
\smallskip

{\em Spectral simplex method}. $A_{k+1}$ is obtained from $A_k$
by changing one of its rows $\ba_i^{(k)}$ so that $(\ba_i^{(k+1)}, \bv_k) > (\ba_i^{(k)}, \bv_k)$,
where $\bv_k$ is a leading eigenvector of $A_k$. 
We fix $i$ and take the element $\ba_i^{(k)} \in \cF_i$ which maximizes the scalar product with the leading eigenvector~$\bv_k$ over all $\ba_i \in \cF_i$. 
The algorithm terminates when no further step is possible, 
i.e., when the matrix $A_k$ is maximal in each row with respect to~$\bv_k$.  

Thus, each iteration makes a one-line correction of the matrix $A_k$
to increase the projection of this line (i.e., row) to the leading eigenvector~$\bv_k$. 
Of course, such a row~$\ba_i^{(k+1)}$ may not be unique. In~\cite{P16}
the smallest index rule was used: each time 
we  take the smallest $i$ for which 
the row $\ba_i^{(k)}$ is not maximal with respect to~$\bv_k$. 
 This strategy provides a very fast convergence to the optimal matrix. 
In~\cite{Akian1} the  convergence was still speeded up by  
 using the {\em pivoting rule}:  
$i$ is the index for which the ratio 
$\frac{(\ba_i^{(k+1)}, \bv_k)}{(\ba_i^{(k)}, \bv_k)}$ is 
maximal. Thus, we always choose the steepest increase of the scalar 
product. One iteration of this method requires 
the exhaustion of all indices~$i=1, \ldots , d$, 
which may be  a disadvantage in case of 
complicated sets~$\cF_i$. 
\smallskip 

{\em The greedy method}. $A_{k+1}$ is obtained by 
replacing all rows of $A_k$ with the maximal 
rows in their sets $\cF_i$ with respect to the leading eigenvector~$\bv_k$. 

So, in contrast to spectral simplex method, we change not only one row
of~$A_k$, but all non-optimal rows, where we can increase   the scalar product $(\ba_i^{(k)}, \bv_k)$. 
If the row~$\ba_i^{(k)}$ is already optimal, we do not change it and set~$\ba_i^{(k+1)} = \ba_i^{(k)}$.  Otherwise we replace it by the row
$\ba_i^{(k+1)} \in \cF_i$ that gives the maximal scalar 
product~$(\ba_i, \bv_k)$ over 
all  $\ba_i \in \cF_i$. 
\bigskip 

{\em The idea behind}. Both the spectral simplex and the greedy methods 
are actually base on the well-known minimax formula for the Perron eigenvalue
(see, for instance,~\cite{BP}): 
$$
\rho(A) \ = \ \min_{i=1, \ldots , d} \, \sup_{\bx > 0}\, \frac{(\ba_i, \bx)}{x_i}\, , 
$$ 
where $A$ is  a non-negative matrix, $\ba_1, \ldots , \ba_d$ are its rows and 
$\bx = (x_1, \ldots , x_d)$ is an arbitrary strictly positive vector. 
This formula makes it possible to use the product structure of a family of matrices
and to actually work with rows separately. See~\cite{Al, NP1} for more details.   
\bigskip 

{\em Realization}. The maximal 
  scalar product  $(\ba_i, \bv_k)$ over all $\ba \in \cF_i$ is found by solving the 
  corresponding convex problem over~${\rm co}(\cF_i)$. If $\cF_i$ is finite, this is done merely by  exhaustion of all elements of~$\cF_i$.
  If $\cF_i$ is a polyhedron, then we find $\ba_i^{(k+1)}$ among its vertices solving an LP problem  by the (usual) simplex method.

  \bigskip

{\em Convergence}. Denote~$\rho(A_k) = \rho_k$. 
  It was shown in~\cite{P16} that both methods are actually relaxations: 
  $\rho_{k+1} \ge \rho_k$ for all~$k$. Moreover, if  $A_k$ is irreducible, 
  then this inequality is strict. Hence, if in all iterations the matrices are irreducible (this is the case, for instance, when all the uncertainty sets~$\cF_i$  are strictly positive), then both methods produce strictly increasing 
  sequence of spectral radii~$\{\rho_k\}_{k\in N}$. Hence, the 
  algorithm never cycles. If all the sets $\cF_i$ are 
  finite or polyhedral, then the set of extreme points of
  the family~$\cF$ is finite, therefore the algorithm  eventually  
  arrives at the maximal spectral radius~$\rho_{\max}$ within finite time. 
  This occurs at the matrix $A_k$ maximal in each row.  
  
  If all~$\cF_i$  are general compact  sets, then it may happen that 
  although a matrix maximal in each row
   exists, it will not be reached  within finite time. 
  In this case the algorithm can be interrupted 
  at an arbitrary step, after which  we apply the a posteriori estimates for $\rho_{\max}$
  defined below. 
  We denote $\bv_k \, = \, \bigl(v_{k, 1}, \ldots , v_{k, d}\bigr)$, so 
   $\ (\ba_i^{(k)}, \bv_k) = \rho_k v_{k, i}$. Then for 
  an arbitrary matrix~$A$ and for its leading eigenvector~$\bv_k$, we define the following values. 
  \begin{equation}\label{eq.si}
  s_i(A)\ = \ 
  \left\{
  \begin{array}{lll}
  \max\limits_{\bb \in \cF_i} \frac{(\bb , \bv_k)}{v_{k, i}}&; & 
  {v_{k, i}} > 0; \\
  +\infty &; & {v_{k, i}} = 0; 
  \end{array}
  \right. \qquad  
  s(A)\ = \ \max_{i = 1, \ldots , d} s_i(A)
  \end{equation}
  Thus, $s_i(A)$ is the maximal 
  ratio between 
  the value $(\ba , \bv_k)$ over all~$\ba \in \cF_i$ and the 
  $i$th component of~$\bv_k$. Similarly, for the minimum: 
  \begin{equation}\label{eq.ti}
  t_i(A)\ = \ 
  \left\{
  \begin{array}{lll}
  \min\limits_{\bb \in \cF_i} \frac{(\bb , \bv_k)}{v_{k, i}}&; & 
  v_{k, i} > 0; \\
  +\infty &; & {v_{k, i}} =  0; 
  \end{array}
  \right. 
  \qquad 
  t(A)\ = \ \min_{i = 1, \ldots , d} t_i(A)
  \end{equation} 
  Then the following obvious estimates for $\rho_{\max}$ are true: 
  \begin{prop}\label{p.20}~\cite{P16}
  For both the spectral simplex method and the greedy method
  (for maximization and for minimization respectively), we have 
  in~$k$th iteration: 
  \begin{equation}\label{eq.apost}
  \rho_k \ \le \ \rho_{\max}\ \le \ s(A_k)\ ; \qquad 
  t(A_k) \ \le \ \rho_{\min}\ \le \ \rho_k\, . 
  \end{equation} 
  \end{prop} 
  In Section~5 we show that at least  for strictly positive matrices,  
  estimates~(\ref{eq.apost}) converge to $\rho_{\max}$ and to $\rho_{\min}$
  respectively with linear rate. This means that for the maximisation problem 
  $|\rho_k - \rho_{\max}| \, \le \, C\, q^k\, $ and for 
  the minimisation problem 
  $\, |\rho_k - \rho_{\min}| \, \le \, C\, q^k$, where $q \in (0,1)$. 
   Moreover, for the greedy method the convergence is locally quadratic, which explains its efficiency in all practical examples. 
 
  However, for sparse matrices the situation is more difficult.
  The equalities~$\rho_{k+1} \ge \rho_k$ are not necessarily strict. 
  The algorithm may stay in the same value of $\rho_k$ for many iterations 
  and even may cycle.  Before  
  addressing this issue, we write a formal procedure for both algorithms. 
  \bigskip 
  
  \begin{center}
\textbf{2.2.  The algorithms}
\end{center}
\bigskip 
  
\noindent \textbf{The spectral simplex method}. 

\smallskip 
  {\em Initialization.} Taking arbitrary $\ba_i \in \cF_i\, , \ i = 1, \ldots , d$,
we form a matrix $A_1 \in \cA$. Denote its   rows by $\ba_1^{(1)}, \ldots , \ba_d^{(1)}$ and its leading eigenvector $\bv_1$.

\smallskip

{\em Main loop. The $k$th iteration}. We have a matrix $A_{k} \in \cA$ composed with rows $\ba_i^{(k)} \in \cF_i,
i = 1, \ldots , d$. Compute its leading eigenvector
$\bv_{k}$ (if it is not unique, take any of them) and for  $i = 1, \ldots d$, find a solution
$\bar \ba_i \in \cF_i$ of the following problem:
\begin{equation}\label{lp-max}
\left\{
\begin{array}{l}
(\ba_i ,\bv_{k}) \ \to \ \max\\
\ba_i \in \cF_i
\end{array}
\right.
\end{equation}
If $\cF_i$ is finite, then this problem is solved by exhaustion; if $\cF_i$ is polyhedral,
then it is solved as an LP problem, and $\bar \ba_i$ is found among its vertices.

If $(\bar \ba_i, \bv_{k}) = (\ba_i^{(k)}, \bv_{k} )$ and $i \le d-1$, then we set $i = i+1$ and solve  problem~(\ref{lp-max})
for the next row. If $i = d$, then the algorithm terminates. In case $\bv_{k} > 0$, the matrix $A_{k}$ is maximal in each
 row, and  $\rho(A_{k})  = \rho_{\max}$, i.e., $A_{k}$ is a solution. If $(\bar \ba_i, \bv_{k}) > (\ba_i^{(k)}, \bv_{k} )$, then we set $\ba^{(k+1)}_i = \bar \ba_i$,
$\ba_j^{(k+1)} = \ba_j^{(k)}$ for all $j\ne i$ and form the corresponding matrix~$A_{k+1}$.
Thus,~$A_{k+1}$ is obtained from $A_k$ by replacing its $i$th row by $\bar \ba_i$.
Then we start $(k+1)$st iteration.
\smallskip

If we obey the pivoting rule, we  solve $d$ problems~(\ref{lp-max})
for all $i = 1, \ldots , d$ and take the row for which the 
value $s_i(A_k)$ defined in~(\ref{eq.si}) is maximal. 
Then~$A_{k+1}$ is obtained from $A_k$ by replacing its $i$th row by $\bar \ba_i$
\smallskip 

{\em Termination}. If $A_k$ is maximal in each row and $\bv_k>0$, then $\rho_{\max} = \rho(A_k)$ is a solution. Otherwise, we get an estimate for $\rho_{\max}$ by
inequality~(\ref{eq.apost}). 
\smallskip 

{\em End.}
\medskip

\begin{remark}\label{r.25}
{\em If each row of $A_k$ makes the maximal scalar product with 
the eigenvector $\bv_k$, but  $\bv_k$ is not strictly positive
(i.e., $A_k$ is ``almost'' maximal in each row), 
then we cannot conclude that $\rho_{\max} = \rho(A_k)$. However, in this case 
the family~$\cF$ is reducible: the space spanned by the vectors with the same 
support as $\bv_k$ is invariant with respect to all matrices from~$\cF$. 
Hence, the problem of maximizing the spectral radius is reduced to two 
problems of smaller dimensions. Therefore, before running the algorithm we 
can  factorize~$\cF$ to 
several irreducible families. The procedure is written in detail 
in~\cite{P16}. However, this way is not very efficient, since 
the case when the final leading eigenvector is not strictly 
positive occurs  rarely. Another way is to factorize the family  after the termination into 
two families as written above and then run the algorithm separately to both of them. 
}
\end{remark}

\medskip 

The simplex method for minimization problem is the same, replacing
maximum by minimum in the problem~(\ref{lp-max}) and omitting the positivity assumption $\bv_{k}>0$ in the
termination of the algorithm.
\smallskip

  \bigskip

\noindent \textbf{The greedy method}
\medskip

{\em Initialization.} Taking arbitrary $\ba_i \in \cF_i\, , \ i = 1, \ldots , d$,
we form a matrix $A_1 \in \cA$ with rows $\ba_1, \ldots , \ba_d$. Take its leading eigenvector $\bv_1$. Denote  $\ba_i^{(1)}\, = \, \ba_i, \ i = 1, \ldots , d$. 

\smallskip

{\em Main loop. The $k$th iteration}. We have a matrix $A_{k} \in \cA$ composed with rows $\ba_i^{(k)} \in \cF_i,
i = 1, \ldots , d$. Compute its leading eigenvector
$\bv_{k}$ (if it is not unique, take any of them) and for  $i = 1, \ldots d$, find a solution
$\bar \ba_i \in \cF_i$ of the problem~(\ref{lp-max}). 
For each $i = 1, \ldots , d$, we do:
\smallskip

{\em If} $(\bar \ba_i, \bv_{k}) = (\ba_i^{(k)}, \bv_{k} )$, then we set 
the $i$th row of $A_{k+1}$ to be equal to that of $A_k$, 
i.e., $\ba_i^{(k+1)} = \ba_i^{(k)}$. 
\smallskip 

 {\em Otherwise}, if  $(\bar \ba_i, \bv_{k}) > (\ba_i^{(k)}, \bv_{k} )$, then we set $\ba^{(k+1)}_i = \bar \ba_i$. 
 
 \noindent We form the corresponding matrix~$A_{k+1}$.
 If the first case took place for all $i$, i.e., 
 $A_{k+1} = A_k$, and if $\bv_k > 0$, then 
 $A_k$ is optimal in each  row. Thus,   
 the answer is $\rho_{\max} = \rho(A_k)$. If $\bv_k$ has some zero components, 
 then the family $\cF$ is reducible. We stop the algorithm and 
 factorize $\cF$ (see Remark~\ref{r.25}).
 Otherwise, go to $(k+1)$st iteration. 
\smallskip 

{\em Termination}. If $A_k$ is maximal in each row, then $\rho_{\max} = \rho(A_k)$ is a solution. Otherwise, we get an estimate for $\rho_{\max}$ by
inequality~(\ref{eq.apost}). 

\smallskip

{\em End.}

\bigskip

The greedy method for minimization problem is the same, replacing
maximum by minimum in the problem~(\ref{lp-max}) and omitting the positivity assumption $\bv_{k}>0$ in the
termination of the algorithm.
\smallskip

As we mentioned, in case of strictly positive uncertainty sets~$\cF_i$ 
both those algorithms work efficiently and, if all~$\cF_i$ are finite or polyhedral, 
they find the optimal matrix within finite time. However, if rows from~$\cF_i$ 
 have some  zero components, the situation may be different. We begin with the 
 corresponding example, then we modify those methods to converge for arbitrary 
 non-negative uncertainty sets. Then in Section~{5} we estimate the rate of convergence.

\bigskip 

\begin{center}
\large{\textbf{3. The cycling and the divergence phenomena}}
\end{center}
\bigskip

Both the spectral simplex method and the greedy method work 
efficiently for strictly positive product families. However, if some 
matrices have zero components, the methods may not work at all. 
For the spectral simplex method, this phenomenon was 
observed in~\cite{P16}. Here is an example for the greedy method. 

The example is for maximising the spectral radius in dimension 3. The starting matrix is strictly positive, 
so the family is irreducible. Moreover, since the staring matrix 
has a simple leading eigenvalue, it follows that 
the spectral simplex method will not cycle for this family~\cite[Theorem 1]{P16}. 
Nevertheless, the greedy method cycles. Moreover, it stabilises on the 
value of the spectral radius $10$, while the real maximum is $12$. 

Thus, this example shows that 
the greedy method may cycle even for an irreducible family and for 
a totally positive starting matrix. And that the value of the objective function on which the algorithm stabilises is far from the optimal value.  

We consider the product family~$\cF$ of $3\times 3$ matrices
defined by the following three uncertainty sets: 
$$
\cF_1 \ = \ \left\{
\begin{array}{lcr}
(1 , & 1, & 1)\\
(0 , & 5, & 10)\\
(0 , & 10, & 5)\\
(12 , & 0, & 0)
\end{array} 
\right\}
\ \qquad 
\cF_2 \ = \ \left\{
\begin{array}{lcr}
(1 , & 1, & 1)\\
(0 , & 10, & 0)
\end{array} 
\right\}\ ; \qquad 
\cF_3 \ = \ \left\{
\begin{array}{lcr}
(1 , & 1, & 3)\\
(0 , & 0, & 10)
\end{array} 
\right\}
$$
So, $\cF_1, \cF_2$, and $\cF_3$ consists of $4, 2$ and $2$ rows 
respectively. Hence, we have in total 
 $4\times 2\times 2 = 16$ matrices. 
Taking the first row in each set we compose the first matrix 
$A_1$ and the algorithm runs as follows: 
$A_1 \, \rightarrow \, A_2 \rightleftarrows \, A_3$, where
$$
\left(
\begin{array}{lcr}
1  & 1 & 1\\
1  & 1 & 1\\
1  & 1 & 3
\end{array} 
\right)
 \xrightarrow
{\bv_1 = (1, 1, 2)^T}
\left(
\begin{array}{lcr}
0  & 5 & 10\\
0  & 10 & 0\\
0  & 0 & 10
\end{array} 
\right)
\begin{array}{c}
\xrightarrow
{\bv_2 = (2, 2, 1)^T}\\
{}\\
\xleftarrow
{\bv_3 = (2, 1, 2)^T}
\end{array}
\left(
\begin{array}{lcr}
0  & 10 & 5\\
0  & 10 & 0\\
0  & 0 & 10
\end{array} 
\right)
$$
The matrix $A_1$ has a unique eigenvector $\bv_1 = (1, 1, 2)^T$. 
Taking the maximal row in each set $\cF_i$ with respect to this eigenvector, 
we compose the  next matrix~$A_2$. It has a multiple 
leading eigenvalue~$\lambda = 10$. Choosing one of its leading eigenvectors 
$v_2 = (2, 2, 1)^T$, we make the next iteration and obtain  
the matrix $A_3$. It also has  a multiple 
leading eigenvalue~$\lambda = 10$. Choosing one of its leading eigenvectors 
$v_2 = (2, 1, 2)^T$, we make the next iteration and come back to 
the matrix $A_2$. The algorithm cycles. The cycling starts on the second iteration, 
the length of the cycle is $2$.  

\textbf{Some conclusions}. The greedy method gets stuck on the cycle 
$A_2 \rightleftarrows \, A_3$  and on the value of the spectral radius~$\rho = 10$. 
However, the maximal spectral radius~$\rho_{\max}$ of the family~$\cF$ 
is not $10$, but $12$. The value~~$\rho_{\max}$   is realized 
for the matrix 
$$
A_4 \ = \ \left(
\begin{array}{lcr}
12  & 0 & 0\\
1  & 1 & 1\\
1  & 1 & 3
\end{array} 
\right)
$$
which is never attained by the algorithm. 

The algorithm never terminates and, when  stopped,  gives a wrong solution: $\rho(A_2) = \rho(A_3) = 10$ instead of $\rho_{\max} = 12$. 
The a posteriori estimate~(\ref{eq.apost}) gives $10\le \rho_{\max}\le 12.5$. 
The error is $25\, \%$, and this is the best the greedy method can achieve for this family.  

Of course, the cycling in this example can easily be overcome by the standard permutation technique. It suffices to add a small positive $\varepsilon$
to each component of vectors in uncertainty sets. However, we have to keep in mind that this is  an illustrative example in a small matrices. Numerical experiments show that for large sparse matrices cycling occurs quite often. In this case the 
perturbation may lead to significant error in the computing of $\rho_{\max}$. 
Standard methods to control the change of the spectral radius are hardly applicable, because we deal with not one matrix but all matrices of the family $\cA$, and  some of them can have ill-conditioned leading eigenvectors. 
 As a result, to guarantee the small change of the spectral radii of 
 matrices from~$\cA$ 
one needs to take too small $\varepsilon$, which leads to the cycling again due to the computation rounding at each step. We discuss this issue in more detail in subsection~4.2.

\bigskip

\begin{center}
\large{\textbf{4. Is the greedy method applicable for sparse matrices?}}
\end{center}
\bigskip

Can the greedy method be be modified to work with sparse matrices?
We need a unified approach for both minimisation and maximisation problems 
which, moreover, would not reduce the efficiency of the method.   To attack this problem we first reveal two main difficulties caused by sparsity. Then we will see how to treat them. We also mention that the same troubles 
occur for the spectral simplex method but they can be  overcome~\cite{P16}. However, those ideas of 
modification are totally  inapplicable for the greedy method. Note that the 
problem of sparsity is different from irreducibility. The latter can easily be 
solved by reduction to several problems of smaller dimensions~\cite{Akian1, NP1}, while the sparsity is more serious.   
\bigskip 

\begin{center}
\textbf{4.1. Two problems: cycling and multiple choice of~$\bv_k$}
\end{center}
\bigskip

In Section~3 we saw that if some vectors from the uncertainty sets 
have zero components, then the greedy method may cycle and, which is still worse, 
may miss the optimal matrix and give a quite rough final estimate of the 
optimal value.  Indeed, numerical  
experiments show that cycling often occurs for sparse matrices, especially 
for minimization problem (for maximizing it is rare but also possible). 
Moreover, the sparser the matrix the more often the cycling occurs.
\smallskip 

The reason of cycling is hidden in multiple leading eigenvectors. 
Recall that an eigenvalue is called {\em simple} if it is a simple (of multiplicity one)
root of the characteristic polynomial. 
The {\em geometrical multiplicity} 
of an eigenvalue is  the dimension of the subspace spanned by 
all  eigenvectors corresponding to that eigenvalue. 
It never  exceeds the {\em algebraic multiplicity} (the multiplicity 
of the corresponding root of the characteristic polynomial). 
We call the subspace spanned by leading eigenvectors {\em Perron subspace}. 
If 
at some iteration we get a matrix $A_k$ with a multiple leading eigenvector, then we have an uncertainty with 
choosing $\bv_k$ from the Perron subspace. An ``unsuccessful''  choice may cause cycling. In the example in Section~3,  
both matrices $A_2$ and $A_3$ have Perron subspaces of dimension~$2$, 
and this is the reason of the trouble. On the other hand,  
if in all iterations the leading eigenvectors~$\bv_k$ are simple,
i.e., their geometric multiplicity are equal to one,  then cycling 
never occurs as the following proposition guarantees. 
\begin{prop}\label{p.10}
If in all iterations of the greedy method,  
the leading eigenvectors~$\bv_k$ are simple, then the method does not cycle
and (in case of finite or polyhedral sets~$\cF_i$) terminates within finite time. 
The same is true for the spectral simplex method.
\end{prop}
This fact follows from Theorem~\ref{th.30}
proved below. Of course, if all matrices $A_k$ are totally positive, 
or at least irreducible, then by the well-known results of the Perron-Frobenius 
theory~\cite[chapter 13, \S 2]{G}, 
all~$\bv_k$ are simple, and  cycling never occurs. 
But how to guarantee the simplicity of leading eigenvectors 
in case of sparse matrices? For the spectral simplex method, 
this problem was solved by the following curious fact: 
\smallskip 

\noindent \textbf{Theorem A.}~\cite{P16} {\em Consider the spectral simplex method in the maximization problem. If the initial matrix~$A_1 \in \cF$ has a simple leading eigenvector, 
then so do all successive matrices~$A_2, A_3, \ldots $ and the method does not cycle.}
\smallskip 

So, in the spectral simplex method, the issue is solved merely by choosing a 
proper initial matrix~$A_1$. In this case there will be no uncertainty in choosing 
the leading eigenvectors $\bv_k$ in all iterations. Moreover, the leading 
eigenvectors will continuously depend on the matrices, hence the 
algorithm will be stable under small perturbations. Note that this holds only for maximization problem! For minimizing the spectral radius,  
the cycling can be avoided as well but in a totally different way~\cite{P16}.  

However, Theorem~A does not hold for the greedy method. This is seen already in the example from Section~3: the totally positive (and hence 
having a simple leading eigenvector) matrix~$A_1$ produces  the matrix~$A_2$ with 
multiple leading eigenvectors. So, a proper choice of the initial matrix~$A_1$ 
is not a guarantee against cycling of the greedy method!

In practice, because of rounding errors, not only multiple 
but ``almost multiple'' leading eigenvectors (when the second largest eigenvalue 
is very close to the first one) may cause cycling or at least essential slowing down 
of the algorithm.   

Nevertheless, for the greedy method, a way to avoid cycling does exist even for  
  sparse matrices. We suggest a strategy  showing its efficiency both 
theoretically and numerically. Even for very sparse matrices,  the algorithm 
works as fast as for positive ones. 

 \bigskip

\begin{center}
\textbf{4.2. Anti-cycling by selection of the leading eigenvectors}
\end{center}
\bigskip

A naive  idea to avoid cycling would be  to slightly 
change all the uncertainty sets~$\cF_i$ making them strictly positive. 
For example, to replace them by  perturbed sets $\cF_{i, \varepsilon} = \cF_i + \varepsilon \be$, where 
$\be = (1, \ldots , 1)\in \re^d$ and $\varepsilon > 0$ is a small constant. So, 
we add a positive number $\varepsilon$ to each entry of vectors from the uncertainty sets. 
All matrices $A \in \cF$
become  totally positive:~$A\, \mapsto \, A_{\varepsilon} = A + \varepsilon E$, where 
$E = \be \be^T$ is the matrix of ones. By Proposition~\ref{p.10}, 
the greedy method for the perturbed family~$\cF_{\varepsilon}$ does not cycle.
However, this perturbation can significantly change the answer of the problem.  
For example, if a matrix $A$ has $d-1$ ones over the main diagonal 
and all other elements are zeros 
($a_{ij} = 1$ if $i - j = 1$ and $a_{ij} = 0$ otherwise), 
then $\rho(A) = 0$, while $\rho(A_{\varepsilon}) > \varepsilon^{1/d}$. 
Even if $\varepsilon$ is very small, say, is $10^{-20}$, 
then for $d=50$ we have  $\rho(A_{\varepsilon}) > 0.4$, while 
$\rho(A) = 0$. 
Hence, we cannot merely replace the family $\cF$ by $\cF_{\varepsilon}$
without risk of a big error. The standard methods to control the change of the spectral radius are not applicable, because we have to control it 
for all matrices of the family $\cA$ simultaneously. Since some of matrices from $\cA$ have ill-conditioned leading eigenvectors (the matrix of eigenvectors have a large 
condition number), etc., this leads to neglectably small $\varepsilon$, which are smaller than the rounding errors.

Our crucial idea is to deal with the original family~$\cF$ but using 
the eigenvectors $\bv_k$ from the perturbed family~$\cF_{\varepsilon}$.  
 \smallskip 
 
\begin{defi}\label{d.20}
The selected leading eigenvector  of a non-negative
matrix $A$ is the limit $\lim\limits_{\varepsilon \to 0}\bv_{\varepsilon}$, 
where $\bv_{\varepsilon}$ is the normalized leading eigenvector of the 
perturbed matrix~$A_{\varepsilon} \, = \, A\, +\, \varepsilon E$. 
\end{defi}
There is a closed formula for the selected eigenvector, which is, 
however, difficult to use in practice~\cite{P19}. Nevertheless, 
the selected eigenvector  can be found explicitly due to the next theorem. 
Before formulating it we make some observations.

If an eigenvector~$\bv$ is simple, then it coincides with the selected eigenvector. 
If~$\bv$ is multiple, then the selected eigenvector is 
a vector from the Perron subspace. Thus, Definition~\ref{d.20} 
provides a way to select 
one vector from the Perron subspace of every non-negative matrix. 
We are going to see that this way is successful for the greedy method.   

Consider the power method for computing 
the leading eigenvectors: $\, \bx_{k+1} \, = \, A\bx_k, \, k \ge 0$. 
In most of practical cases it converges to a leading eigenvector~$\bv$ linearly 
with the rate $O \bigl( \bigl| \frac{\lambda_2}{\lambda_1}\bigr|^k\bigr)$, 
where $\lambda_1, \lambda_2$ are the first and the second largest by modulus 
eigenvalues of~$A$ or $O\bigl( \frac{1}{k} \bigr)$, if the leading eigenvalue 
has non-trivial (i.e., of size bigger than one) Jordan blocks. 
In a rare case  when there are several largest by modulus eigenvalues, 
then the ``averaged'' sequence 
$\frac1r\, \sum_{j=0}^{r-1} x_{k+j}$ converges to $\bv$, 
where $r$ is the   {\em imprimitivity index} (see~\cite[chapter 13]{G}). 
In all cases we will say that the power method converges and, for the sake of simplicity, assume that there is only one largest by modulus eigenvalue, maybe multiple, i.e., that $r=1$.   

\begin{theorem}\label{th.40}
The selected leading eigenvector of a non-negative
matrix $A$  is proportional to the limit of the power method 
applied to the matrix~$A$ with the initial vector~$\bx_0 = \be$. 
\end{theorem}
\noindent {\tt Proof.} It may be assumed that $\rho(A) = 1$. 
The spectral radius of the matrix $A_{\varepsilon}$
strictly increases in $\varepsilon$, so  
$\rho(A_{\varepsilon}) =  1 + \delta$ for some $\delta > 0$. We have 
$$
\bigl( A\, + \, \varepsilon \, \be \, \be^T \bigr)\, \bv_{\varepsilon} \ = \ 
(1+\delta)\bv_{\varepsilon}\, .  
$$
Denoting by $S_{\varepsilon} = (\be, \bv_{\varepsilon})$ the 
sum of entries of the vector~$\bv_{\varepsilon}$, we obtain 
$$
A\bv_{\varepsilon}\ + \ \varepsilon \, S_{\varepsilon}\, {\be}   \ = \ 
(1+\delta)\, \bv_{\varepsilon}\, , 
$$ 
and hence 
$$
\left( I \, - \, \frac{1}{(1+\delta)}\, A   \right)\,  \bv_{\varepsilon}
\quad = \quad \frac{\varepsilon \, S_{\varepsilon}}{1+\delta} \, \be \, .
$$ 
Since $\rho(A) = 1$, we have 
$\rho\bigl(\frac{1}{(1+\delta)}\, A \bigr)\, = \, \frac{1}{1+\delta} \, < \, 1$.
Therefore,  
we can apply the power series expansion: 
$$
\left( I \, - \, \frac{1}{(1+\delta)}\, A   \right)^{-1} \ = \ 
I \, + \, (1+\delta)^{-1}A \, + \,  (1+\delta)^{-2}A^2\, + \, \ldots  
$$ 
 and obtain 
$$
\bv_{\varepsilon} \quad  = \quad  
\frac{1+\delta}{\varepsilon \, S_{\varepsilon}}\ \sum_{k=0}^{\infty} 
\ (1+\delta)^{-k}A^k {\be}\, .
$$
Assume now that the power method for the matrix~$A$ and for $\bx_0 = \be$ 
converges to some vector~$\tilde \bv$. In case $r=1$ ($r$ is the imprimitivity index), this means that $A^k\be\to \tilde \bv$
as $k \to \infty$. Then 
  direction of the 
vector  $\sum_{k=0}^{\infty} (1+\delta)^{-k}A^k {\be}$ converges 
as $\delta \to 0$ to the direction of the vector~$\tilde \bv$. 
Since $\delta \to 0$ as $\varepsilon \to 0$, the theorem follows. 
If $r>1$, then $\frac1r\, \sum_{j=0}^{r-1} A^{k+j}\be\to \tilde \bv$. 
Arguing as above we  conclude again that the direction of the 
vector  $\sum_{k=0}^{\infty} (1+\delta)^{-k}A^k {\be}$ converges 
as $\delta \to 0$ to the direction of the vector~$\tilde \bv$.
This completes the proof.

{\hfill $\Box$}
\smallskip

\begin{defi}\label{d.40}
The greedy method with the selected leading eigenvectors~$\bv_k$ 
in all iterations is called selective greedy method.
\end{defi}

Now we arrive at the main result of this section. We show 
that using selected eigenvectors~$\bv_k$ avoids cycling. 

\begin{theorem}\label{th.30}
The selective greedy method does not cycle.
\end{theorem}
\noindent {\tt Proof.} Consider the case of finding maximum, 
the case of minimum is proved in the same way.
By the $k$th iteration of the algorithm we have a matrix $A_k$ and its 
selected leading  eigenvector~$\bv_k$. Assume the algorithm cycles: 
$A_{1} = A_{n+1}$, and let $n\ge 2$ be the minimal length of the cycle
(for $n=1$, the equality $A_{1} = A_{2}$ means that $A_1$ is the optimal matrix, 
so this is not cycling).  
 Consider the chain of 
perturbed matrices $A_{1, \varepsilon}\, \to \, \cdots 
\to A_{n+1, \varepsilon}$, where  $A_{k, \varepsilon}\, =
\, A_{k} \, +\, \varepsilon \, E$
and $\, \varepsilon >0$ is a small number. For the matrix $A_{k}$, in each 
row $\ba^{(k)}_i$, we have one of two cases: 
\smallskip 

1)  if $\bigl( \ba^{(k)}_i, \bv_k\bigr) \, = \, 
\max\limits_{\ba_i \in \cF_i}(\ba_i,
\bv_k)$, then this row is not changed in this iteration: 
$\ba^{(k+1)}_i \, = \, \ba^{(k)}_i$. In this case
$\ba^{(k+1)}_i\, +\, \varepsilon \be  \, = \, \ba^{(k)}_i + \varepsilon \be$. 
\smallskip 

2) if $\bigl(\ba^{(k)}_i, \bv_k\bigr) \, < \, \max\limits_{\ba_i \in \cF_i}(\ba_i,
\bv_k)$, then the row $\ba^{(k)}_i$ is not optimal, 
and $\bigl(\ba^{(k)}_i, \bv_k\bigr) \, < \, \bigl( \ba^{(k+1)}_i, \bv_k
\bigr)$.  
Hence the same inequality is true for the perturbed 
matrices $A_{k, \varepsilon}, A_{k+1, \varepsilon}$
and for the eigenvector $\bv_{k, \varepsilon}$
of $A_{k, \varepsilon}$, whenever $\varepsilon$ is 
small enough. 
\smallskip 

Thus, in both cases, each row of 
$A_{k+1, \varepsilon}$ makes a bigger or  equal
scalar product with $\bv_{k, \varepsilon}$ than 
the corresponding row of $A_{k, \varepsilon}$. 
Moreover, at least in one row this scalar product is strictly 
bigger, otherwise $A_{k} = A_{k+1}$, which contradicts 
the minimality of the cycle. This, in view of strict positivity of 
$A_{k, \varepsilon}$, implies that 
$\rho_{k+1, \varepsilon} \, > \, \rho_{k, \varepsilon}$~\cite[lemma 2]{P16}. We see that in the chain $A_{1, \varepsilon}\to \cdots 
\to A_{n+1, \varepsilon}$ the spectral radius strictly  increases.   
Therefore $A_{1, \varepsilon}\ne 
 A_{n+1, \varepsilon}$, and hence $A_1\ne 
 A_{n+1}$. The contradiction competes the proof.

{\hfill $\Box$}
\smallskip

\begin{ex}\label{ex.10}
{\em We apply the selective greedy method to the family from Section~3, 
on which the (usual) greedy method cycles and arrives to a wrong solution. 
Now  we have: 
$$
A_1 \ \xrightarrow{\ } \ A_2 \ \xrightarrow{\bv_2 = (3, 2, 2)^T}\ 
\left(
\begin{array}{lcr}
12  & 0 & 0\\
0  & 10 & 0\\
0  & 0 & 10
\end{array} 
\right)
 \ \xrightarrow
{\bv_3 = (1, 0, 0)^T}\ 
A_4 \ = \ \left(
\begin{array}{lcr}
12  & 0 & 0\\
1  & 1 & 1\\
1  & 1 & 3
\end{array} 
\right)\ 
\xrightarrow
{\bv_4 = (49, 5, 6)^T} \ 
A_4
$$
Thus, the selective greedy method arrives at the optimal matrix~$A_4$
at the third iteration. This matrix is maximal in each row with respect to 
its leading eigenvector~$\bv_4$, as the last iteration demonstrates. 

The selective greedy method repeats the first iteration of the 
(usual) greedy method, but in the second and in the third iteration it chooses different leading eigenvectors (``selected eigenvectors'')
by which it does not cycle and goes directly to the optimal solution. 
}

\end{ex}
\bigskip 

\newpage 

\begin{center}
\textbf{4.3. Realization of the selective greedy method}
\end{center}
\bigskip

Thus, the greedy method does not cycle, provided  it uses selected leading 
eigenvectors~$\bv_k$ in all iterations. 
If the leading eigenvector is unique, then 
it coincides with the selected  leading eigenvector, and 
there is no problem. The difficulty occurs if the leading eigenvector is multiple.
The following crucial fact is a straightforward consequence of 
Theorem~\ref{th.40}. 

\begin{cor}\label{c.20}
The power method $\bx_{k+1} = A\bx_k, \, k \ge 0,$ applied to the  initial vector $\bx_0 = \be$ converges to the selected leading eigenvector. 
\end{cor}

Therefore, to realize the selective greedy method  we compute all the
leading eigenvectors~$\bv_k$ by the power method starting always with 
the same vector $\bx_0 = \be$ (the vector of ones). After a sufficient number of iterations we come close to the limit, which is, by Corollary~\ref{c.20}, the 
selected eigenvector (up to normalization). So, actually we perform the greedy method with 
approximations for the selected eigenvectors, which are, in view 
of Theorem~\ref{th.30}, the leading eigenvectors of 
perturbed matrices $A_{\varepsilon}$ for some small~$\varepsilon$.   

Thus, to avoid cycling we can compute all eigenvectors $\bv_k$
by the power method starting with the same initial vector~$\bx_0 = \be$. 
Because of computer approximations, actually we deal with 
totally positive matrices of the family~$\cF_{\varepsilon}$ for some 
small~$\varepsilon$.  

\begin{remark}\label{r.123}
{\em The usual method for the non-symmetric eigenvalue/eigenvector problem 
is the QR factorisation. The case of non-negative matrices is not an exception. 
However, for large and for sparse non-negative matrices, in many cases the Perron eigenvalue problem 
is solved more efficiently by the power method. We use power method to avoid 
cycling of the greedy algorithm. Namely, the power method allows us to find
selected leading eigenvector, see~\cite{P19}.   
}
\end{remark}
\bigskip

\begin{center}
\textbf{4.4. The procedure of the selective greedy method}
\end{center}
\bigskip

We precisely follow the algorithm of the greedy method presented in 
subsection~2.2, and in each iteration we compute the corresponding  leading eigenvectors~$\bv_k$ by the power method starting with the vector of ones~$\be$. 
The precision parameter $\varepsilon$ of the power method is chosen 
in advance and is the same in all iterations.  By Theorem~\ref{th.30}, 
if $\varepsilon > 0$ is small enough, then the greedy method does not cycle and 
finds the solution within finite time. In practice, if the algorithm cycles for 
a given $\varepsilon$, then we reduce it taking, say, $\varepsilon/10$ and restart the algorithm. 

It is known that the power method may suffer for non-primitive  matrices. In this case, 
the eigenvalue with the biggest modulus may not be unique, which can cause slow down or even divergence of the algorithm. In particular, this explains possible difficulties 
of the power method for sparse matrices. See, for example,\cite{NN}. Here to avoid any trouble with the power method for sparse matrices, we can use the following well-known fact: if $A\ge 0$ is a matrix, then the matrix
$A+I$ has a unique (maybe multiple) eigenvalue equal by modulus to 
$\rho(A)+1$. This implies that the power method applied for the matrix $A+I$
always converges to its leading eigenvector, i.e., $(A+I)^k\be \to \bv$ as $k\to \infty$. 
Of course, $A+I$ has the same leading eigenvector as $A$. 
Hence, we can replace each  uncertainty set~$\cF_i$ 
by the set $\be_i + \cF_i$, where $\be_i$ is the $i$th canonical 
basis vector. After this each matrix $A_k$ is replaced by $A_k+I$, and 
we will not have trouble with the convergence of the power method. 

\bigskip

\begin{center}
\large{\textbf{5. The rate of convergence}}
\end{center}
\bigskip

In this section we address two issues:  explaining  the fast 
convergence of the greedy method and realising why the spectral simplex method 
converges slower. We restrict ourselves 
to the maximization problems, where we find the maximal spectral radius~$\rho_{\max}$. 
For minimization problems, the results and the proofs are the same.   

First of all, let us remark that due to possible cycling, no efficient 
estimates for the rate of convergence exist. Indeed, in Section~3 
we see that the greedy method may not converge to the solution at all. 
That is why we estimate the rate of convergence only under some favorable 
assumptions (positivity of the limit matrix, positivity of its leading eigenvector, etc.) Actually, they are  not very restrictive since  the selective greedy method has the same convergence rate as the greedy method for a strictly positive family~$\cF_{\varepsilon}$.

We are going to show  that under those assumptions, both methods 
have global linear convergence and the greedy method has moreover a local quadratic convergence. In both cases we estimate the parameters of linear and quadratic convergence.

\bigskip

\begin{center}
\textbf{5.1. Global linear convergence}
\end{center}
\bigskip

The following result is formulated for both the greedy and the spectral simplex methods. The same holds for the method with the pivoting rule. We denote by $A_k$ the matrix obtained by the $k$th iteration, 
by $\rho_k$ its spectral radius and by $\bu_k, \bv_k$ its  left and right 
leading eigenvectors (if there are multiple ones, we take any of them).
  The $j$th components of the vectors $\bu_k, \bv_k$ are 
   $u_{k, j}, \, v_{k, j}$ respectively.   
We also denote $\rho_{\max} = \max_{A \in \cF}\rho(A)$. As usual, 
 $\{\be_j\}_{j=1}^d$ is the canonical basis in $\re^d$ and $\be$ is the vector 
 of ones.  
\begin{theorem}\label{th.50}
In both the greedy method and the spectral simplex method for maximizing the spectral radius, in the $k$th   iteration,  we have  
\begin{equation}\label{eq.estimate}
\rho_{k+1} \quad   \ge  \quad  \rho_k \ + \ (\rho_{\max} - \rho_k)\, 
\max\limits_{j = 1,\ldots,d}\frac{u_{k+1, j} \, v_{k, j} }{(\bu_{k+1}, \bv_k)}\ ,
\qquad k \in \n,     
\end{equation}
provided $\bv_k > 0$. 
\end{theorem}
\noindent {\tt Proof.} Consider the value $s = s(A)$ defined in~(\ref{eq.si}). 
Let the maximum in formula~(\ref{eq.si}) be attained in the
$m$th row.  
Then for every matrix $A\in \cF$, we have $A\bv_k \, \le \, s\, \bv_k$, and therefore $\rho(A) \le s$. Thus, $\rho_{\max} \le s$. 
On the other hand,  in both 
greedy and simplex methods, the $m$th component of the vector $A_{k+1}\bv_k$
is equal to the $m$th component of the vector $\, s\, \bv_k$. 
In all other components, we have 
$A_{k+1}\bv_k \ge A_{k}\bv_k = \rho_k\, \bv_k$. 
Therefore,  
$$
A_{k+1}\bv_k \ \ge \ \rho_k \bv_k \ + \ (s-\rho_k) v_{k, j}\be_j\, . 
$$
Multiplying this inequality by the vector $\bu_{k+1}$ from the left, we obtain 
$$
\bu_{k+1}^TA_{k+1}\bv_k \ \ge \ \rho_k (\bu_{k+1}, \bv_k) \ + \ (s-\rho_k)\, 
v_{k, j}\, (\bu_{k+1}, \be_j)\, . 
$$
Since $\, \bu_{k+1}^TA_{k+1} \, = \, \rho_{k+1}\bu_{k+1}^T\, $ and 
$(\bu_{k+1}, \be_j) \, = \, u_{k+1, j}$, we have 
$$
\rho_{k+1}(\bu_{k+1}, \bv_k) \ \ge \ \rho_k (\bu_{k+1}, \bv_k) \ + \ 
(s-\rho_k) v_{k, j} u_{k+1, j}\, . 
$$
Dividing by $(\bu_{k+1}, \bv_k)$ and taking into account that $s\ge \rho_{\max}$, we arrive at~(\ref{eq.estimate}). 

{\hfill $\Box$}
\smallskip

Rewriting~(\ref{eq.estimate}) in the form 
$$
\rho_{\max}\ - \ \rho_{k+1} \quad \le \quad 
   \ (\rho_{\max}\ - \ \rho_{k})\, \left(  
1 \, - \, \frac{u_{k+1, j} \, v_{k, j}}{(\bu_{k+1}, \bv_k)}\right)\, ,   
$$
we see that the $k$th iteration reduces the distance 
to the maximal spectral radius in at least 
$\left(  1 \, - \, \frac{u_{k+1,j}\, v_{k,j}}{(\bu_{k+1}, \bv_k)}\right)$
times. Thus, we have established  
\begin{cor}\label{c.10}
Each  iteration of the greedy method and of the spectral simplex method
reduces the distance to the optimal spectral radius 
in at least $\left(  1 \, - \, \frac{u_{k+1, j}\, v_{k,j}}{(\bu_{k+1}, \bv_k)}\right)$
times. 
\end{cor}
If  $A_k > 0$, then the contraction coefficient from Corollary~\ref{c.10}
can be roughly estimated from above. Let $m$ and $M$ be the smallest and the biggest 
entries respectively of all vectors from the uncertainty sets. 
For each $A\in \cF$,  we denote by $\bv$ its leading eigenvector and by $\rho$
its spectral radius. Let also $S$ be the sum of entries of~$\bv$. 
Then for each $i$, we have $v_i = \rho^{-1}(A\bv, \be_i) \, \le \, \rho^{-1}MS$. 
Similarly, $v_i \ge \rho^{-1}mS$. Hence, for every matrix from $\cF$,  
the ratio between the smallest and of the biggest entries of  its leading eigenvector 
is at least $m/M$. The same is true for the left eigenvector. 
Therefore, 
$$
\frac{u_{k+1, j}v_{k, j}}{(\bu_{k+1}, \bv_k)} \quad \ge \quad \frac{m^2}{m^2 + (d-1)M^2}.
$$
We have arrived at the following theorem on the global linear convergence  
\begin{theorem}\label{th.60}
If all uncertainty sets are strictly positive, 
then the rate of convergence of both the greedy method and of the spectral simplex method is at least linear with the factor 
\begin{equation}\label{eq.q}
q \quad \le \quad 1 \ - \  
\frac{m^2}{m^2 + (d-1)M^2}\, , 
\end{equation}
 where $m$ and $M$ are 
the smallest and the biggest 
entries respectively of vectors from the uncertainty sets.
\end{theorem}
Thus, 
in each  iteration of the greedy method and of the spectral simplex method, we have  
\begin{equation}\label{eq.linear}
\rho_{\max}\,- \, \rho_{k+1} \ \le  \ 
q^{k}\, (\rho_{\max} \, - \, \rho_0),   
\end{equation}

\bigskip

\newpage 

\begin{center}
\textbf{5.2. Local quadratic convergence}
\end{center}
\bigskip

Now we are going to see that the greedy method 
have actually a quadratic rate of convergence in a 
small neighbourhood of the optimal matrix. This explains its fast convergence. We establish this result under  the following two assumptions: 1) the optimal matrix $\bar A$, for which $\rho(\bar A) = \rho_{\max}$, is
strictly positive; 2) in a neighbourhood of $\bar A$, 
 the uncertainty sets $\cF_i$ are bounded by $C^2$ smooth strictly 
convex surfaces. The latter is, of course, not satisfied for polyhedral sets. 
Nevertheless, the quadratic convergence for sets with a smooth boundary explains 
the fast convergence for finite or for polyhedral sets as well.      

The  local quadratic convergence means  that 
there are constants $B > 0$ and $\varepsilon \in \bigl(0, \frac{1}{B} \bigr)$
for which the following holds: 
if $\|A_k-\bar A\| \le \varepsilon$, then $\|A_{k+1}-\bar A\| \, \le \, 
B\, \|A_k-\bar A\|^2$. In this case, the algorithm converges whenever 
$\|A_1-\bar A\| \le  \varepsilon\, $, in which case  for every 
iteration~$k$, we have 
$$
\|A_k-\bar A\| \ \le \ \frac{1}{B}\, \left( B\, \|A_1-\bar A\|\right)^{2^k-1}
$$ 
(see, for instance,~\cite{A}). 

As usual, we use Euclidean norm~$\|\bx\| = \sqrt{(\bx, \bx)}$. We also need 
the $L_{\infty}$-norm $\|\bx\|_{\infty} = \max_{i=1, \ldots , d}|x_i|$. 
We denote by $I$ the identity $d\times d$ matrix; $A \preceq B$ means that 
the matrix $B-A$ is positive semidefinite.  
A set $\cM$ is {\em strongly convex} if its boundary does not contain segments.  
\begin{theorem}\label{th.10}
Suppose the optimal matrix $\bar A$ is strictly positive
and all uncertainty sets $\cF_i$ are strongly convex and 
$C^2$ smooth in a neighbourhood of~$\bar A$; 
then the greedy algorithm has a local quadratic convergence to $\bar A$. 
\end{theorem}
In the proof we use the following auxiliary fact
\begin{lemma}\label{l.10}
If $A$ is a strictly positive matrix
with $\rho(A) = 1$ and
with the leading eigenvector $\bv$, then for every 
$\varepsilon > 0$ and for every 
 vector $\bv'\ge 0, \, \|\bv'\| = 1$, 
such that $\|A\bv' - \bv'\|_{\infty} \le \varepsilon$,  we have 
$\|\bv - \bv'\|\, \le \, C(A) \, \varepsilon$, where 
$C(A) \, = \, \frac{1}{m} \, \bigl(1 + (d-1)\frac{M^2}{m^2}\bigr)^{1/2}$, 
$m$ and $M$ are minimal and maximal entry of $A$ respectively. 
\end{lemma}
{\tt Proof.} Let $\bv = (v_1, \ldots , v_d)$. 
Denote $\bv' = (v_1s_1, \ldots , v_ds_d)$, where $s_i$ are some non-negative numbers. Choose one for which the value $|s_i - 1|$ is maximal. Let it be $s_1$ and let $s_1 - 1 >0$
(the opposite case is considered in the same way). 
Since $\|\bv'\| = \|\bv\|$ and $s_1> 1$, it follows that there exists 
$q$ for which $s_q< 1$. 
Since $A\bv = \bv$, we have $\sum_{j=1}^d a_{1j}v_j = v_1$. 
Therefore, 
\begin{equation}\label{eq.est10}
(\bv' - A\bv')_{1} \  = \ s_1v_1 \ - \  \sum_{j=1}^d a_{1j}v_js_j \ = \ 
\sum_{j=1}^d a_{1j}v_js_1 \ - \  \sum_{j=1}^d a_{1j}v_js_j \ = \ 
\sum_{j=1}^d a_{1j}v_j(s_1 - s_j)
\end{equation}
The last sum is bigger than or equal to one term 
$a_{1q}v_m(s_1 - s_q) \, \ge \,  a_{1q}v_q(s_1 - 1) \, \ge \, m (v_1' - v_1)$.   
Note that for all other $i$, we have $|s_i - 1| \le |s_1 - 1|$, hence 
$|v_i' - v_i| \le \frac{v_i}{v_1}|v_1' - v_1| \le \frac{M}{m}|v_1' - v_1|$. 
Therefore, $\|\bv' - \bv\|\, \le \, \bigl(1 + (d-1)\frac{M^2}{m^2}\bigr)^{1/2}|v_1' - v_1|$. Substituting to~(\ref{eq.est10}), we get 
$$
(\bv' - A\bv')_{1} \ \ge \ m |v_1' - v_1| \ \ge \ m \left(1 + (d-1)\frac{M^2}{m^2}\right)^{-1/2}\|\bv' - \bv\|
$$
Hence, $\|A\bv' - \bv'\|_{\infty} \ge \frac{1}{C(A)} \|\bv' - \bv\|$, which completes the proof. 

{\hfill $\Box$}
\medskip 

{\tt Proof of Theorem~\ref{th.10}.} 
Without loss of generality it can be assumed that $\rho(\bar A) = 1$. 
It is known that any strongly convex smooth surface
$\Gamma$ can be defined by  a smooth convex function $f$
so that $\bx \in \Gamma \Leftrightarrow f(\bx)=0$ and $\|f'(\bx)\| = 1$
for all $\bx \in \Gamma$. For example, one can set $f(\bx)$ equal to  the distance 
from $\bx$ to $\Gamma$ for $\bx$ outside $\Gamma$ or in some neighbourhood inside it. 
For each $i = 1, \ldots , d$, 
we set $\Gamma_i = \partial \cF_i$ and denote by $f_i$ such a function 
 that 
defines the surface $\Gamma_i$.  
Fix an index $i = 1, \ldots , d$ and 
denote by $\ba^{(k)}, \ba^{(k+1)},$ and $\bar \ba$  the~$i$th rows of 
the matrices $A_k, A_{k+1},$ and $\bar A$ respectively
(so, we omit  the subscript~$i$ to simplify the notation).
Since $\bar \ba \, = \, {\rm argmax}_{\ba \in \Gamma_i}\, (\ba , \bar \bv)$, 
it follows that $f_i'(\bar \ba) = \bar \bv$.  
Denote by $L$ and $\ell$
the lower and upper bounds of the quadratic form $f''(\ba)$ in some  
neighbourhood of $\bar \ba$, i.e., 
$\, \ell \, I \preceq f''(\ba) \preceq L\, I$ at all points 
 $\ba$ such that 
$\|\ba - \bar \ba\| < \varepsilon$.  
Writing the Tailor expansion of the function $f_i$ 
at the point~$\bar \ba$, we have for small $\bh\in \re^d$: 
$$
f_i(\bar \ba + \bh)\ = \ f_i(\bar \ba )\, + \, 
\bigl( f_i'(\bar \ba) , \bh\bigr) \, + \, 
r(\bh)\|\bh\|^2, 
$$ 
where the remainder $r(\bh) \le L$. Substituting 
$\bh = \ba^{(k)} - \bar \ba$ and taking into account that 
$f_i(\ba^{(k)})\, = \, f_i(\bar \ba ) \, = \, 0$, because 
both $\ba^{(k)}$ and $\bar \ba$ belong to $\Gamma$, we obtain  
$\bigl( f_i'(\bar \ba) , \bh\bigr) \, + \, r(\bh)\|\bh\|^2\, = \, 0$. 
Hence $\, |(\bar \bv , \bh)|\, \le \, L \|\bh\|^2$. 
Thus, 
$$
|(\bar \bv , \ba^{(k)}) \, - \, (\bar \bv , \bar \ba)| \ \le \ L\,  \|\bh\|^2\, .   
$$
This holds for each row of $A_k$, therefore
$\|(A_k - I)\bar \bv\|_{\infty}\, \le \, L \|\bh\|^2$, and hence 
$\|\bv_k - \bar \bv\| \, \le \, C(A_k)\, L\, \|\bh\|^2$, where 
the value $C(A_k)$ is defined in Lemma~\ref{l.10}. 

Further, $\ba^{(k+1)} \, = \, {\rm argmax}_{\ba \in \Gamma_i}\, (\ba ,  \bv_k)$, hence 
$f_i'(\ba^{(k+1)}) \,  = \, \bv_k$.  Consequently, 
$$
\|\bv_k - \bar \bv\| \ = \ |f_i'(\ba^{(k+1)}) - f_i'(\bar \ba)|\ \ge \ 
\ \ell \, \|\ba^{(k+1)} - \bar \ba\|. 
$$
Therefore, 
$\|\ba^{(k+1)} - \bar \ba\| \, \le \, \frac{C(A_k)L}{\ell} \|\bh\|^2\, = \, 
\frac{C(A_k)L}{\ell} \|\ba^{(k+1)} - \bar \ba\|^2$. This holds for 
each row of the matrices, hence the theorem follows.

{\hfill $\Box$}
\medskip

In the proof of Theorem~\ref{th.10} we estimated  the parameter $B$ of the  quadratic convergence in terms of the parameters $\ell$ and $L$ of the 
functions $f_i''$.    
 This estimate, however, is difficult to evaluate, because the 
 functions $f_i$ are a priori unknown. 
To estimate  $B$ effectively, we need  another 
idea based on geometrical properties of the uncertainty sets $\cF_i$. 
Below we obtain such an estimate in terms of radii of curvature of 
the boundary of $\cF_i$. 

Assume that at each point  of intersection of the surface $\Gamma_i = \partial \cF_i$ 
and of  the corresponding $\varepsilon$-neighbourhood of the point $\bar \ba_i$, 
the maximal radius of curvature of two-dimensional cross-section passing through the 
normal vector to the surface 
does not exceed  $R$ and the  minimal radius of curvature is 
at least $r > 0$. Denote also by $C$ the maximal value of 
$C(A)$ over the corresponding neighbourhood of the matrix 
$\bar A$, where $C(A)$ is defined in Lemma~\ref{l.10}. 

\begin{theorem}\label{th.20}
For the greedy method, we have $\, B \, \le \, \frac{RC}{2r \rho(\bar A)}$. 
\end{theorem}
{\tt Proof.} As in the proof of Theorem~\ref{th.10}, we assume that 
$\rho(\bar A) = 1$ and denote by 
$\ba^{(k)}, \ba^{(k+1)},$ an $\bar \ba$ the~$i$th rows of 
the matrices $A_k, A_{k+1},$ and $\bar A$ respectively.
Since $\bar \ba \, = \, {\rm argmax}_{\ba \in \Gamma_i}\, (\ba , \bar \bv)$, 
it follows that $\bar \bv$ is the outer unit normal vector to~$\Gamma_i$ at the 
point~$\bar \ba$.   Denote $\bh = \ba^{(k)} - \bar \ba$. 
Draw a Euclidean sphere~$\cS_i$ of radius $r$ tangent to $\Gamma_i$ from inside at 
the point $\bar \ba$. Take a point  $\bb^{(k)}$  on $\cS_i$ such that
$\|\bb^{(k)} - \bar \ba\| = \|\bh\|$. Since the maximal radius of curvature of~$\Gamma_i$
is at least $R$, this part of the  sphere is inside~$\cF_i$, hence 
the vector $\ba^{(k)} - \bar \ba$ forms a smaller angle with 
the normal vector $\bar \bv$ than the vector $\bb^{(k)} - \bar \ba$. 
The vectors $\ba^{(k)} - \bar \ba$ and $\bb^{(k)} - \bar \ba$ 
 are of the same length, therefore 
$$
\bigl( \ba^{(k)} - \bar \ba \, , \, \bar \bv\bigr) \ \ge \ 
\bigl( \bb^{(k)} - \bar \ba \, , \, \bar \bv\bigr)\ = \ - \frac{\|\bh\|^2}{2r}
$$
On the other hand,  $\bigl( \ba^{(k)} - \bar \ba  ,  \bar \bv\bigr) \le 0$. Therefore, 
$$
\Bigl|\, \bigl( \ba^{(k)}  \, , \, \bar \bv\bigr) \ - \ 
\bigl(\bar  \ba \, , \, \bar \bv\bigr)\, \Bigl|  \ \le \ 
  \frac{\|\bh\|^2}{2r}
$$
This inequality holds for all rows of the matrix $A_k$, therefore
$$
\|(A_k - I)\bar \bv\|_{\infty}\ = \ 
\|(A_k - \bar A)\bar \bv\|_{\infty}\, \le \, \frac{\|\bh\|^2}{2r}, 
$$ 
and hence, by Lemma~\ref{l.10},  
\begin{equation}\label{eq.vk2}
\|\bv_k - \bar \bv\| \quad \le \quad  \frac{C }{2r}\ \|\bh\|^2\, .
\end{equation}
Further, $\ba^{(k+1)} \, = \, {\rm argmax}_{\ba \in \Gamma_i}\, (\ba ,  \bv_k)$, hence 
$\bv_k$ is a normal vector to $\Gamma_i$ at the point~$\ba^{(k+1)}$.   
Consequently  
$$
R\, \|\bv_k - \bar \bv\| \ \ge \ 
\  \|\ba^{(k+1)} - \bar \ba\|. 
$$
Combining this inequality with~(\ref{eq.vk2}) we obtain 
$\|\ba^{(k+1)} - \bar \ba\| \, \le \, \frac{C R}{2r} \, \|\bh\|^2$, 
which concludes the proof.

{\hfill $\Box$}
\medskip 

\bigskip

\begin{center}
\large{\textbf{6. Numerical results}}
\end{center}
\bigskip

We now demonstrate and analyse  numerical results of the selective greedy method.
Several types of problems are considered:  maximizing and minimizing the 
spectral radius, finite and polyhedral uncertainty sets, positive and sparse matrices. 
In all tables $d$ denotes the dimension of the matrix and $N$ denotes the number of 
elements in each uncertainty set~$\cF_i, \, i = 1, \ldots , d$. 
For the sake of simplicity, in each experiment all uncertainty sets~$\cF_i$ 
are of the same cardinality. For each pair $(d, N)$ wee made ten experiments and put 
the arithmetic mean in the table. 
 All computations are done on a standard laptop 
(Dell XPS 13 with Intel Core i7-6500U CPU @ 2.50GHz and 8GB RAM.) 
The algorithms is coded in Python.
\medskip 

\begin{center}
\textbf{6.1. Product families with finite uncertainty sets}
\end{center}
\medskip 

We first test the selective greedy method on positive product families, and then 
 on non-negative product families with  density parameters $\gamma_i\in(0,1)$
  (the percentage of nonzero entries of the vectors belonging to the uncertainty set $\mathcal{F}_i$.)

The results of tests on positive product families are given in Tables 2 and 3, for maximization and minimization respectively. The numbers shown represent the average number of iterations performed by the algorithm to find the optimal matrix.

\begin{table}[h!]
\begin{center}
\begin{tabular}{c|c c c}
 $d\setminus N$ & 50  & 100 & 250\\
\hline
25 & 3.2 & 3 & 3.2\\
100 & 3 & 3.3 & 3\\
500 & 3.1 & 3.1 & 3.2\\
2000 & 3 & 3 & 3.1
\end{tabular}
\caption*{{\footnotesize Table 2: Average number of iterations for maximization, for positive families}}
\end{center}
\end{table}

\begin{table}[h!]
\begin{center}
\begin{tabular}{c|c c c}
 $d\setminus N$ & 50  & 100 & 250\\
\hline
25 & 3.2 & 3.3 & 3.2\\
100 & 3 & 3.1 & 3.2\\
500 & 3.1 & 3.1 & 3\\
2000 & 3 & 3.2 & 3.1
\end{tabular}
\caption*{{\footnotesize Table 3: Average number of iterations for minimization, for positive families}}
\end{center}
\end{table}

In Tables 4 and 5 we report the behaviour of selective greedy algorithm as dimension $d$ and the size of the uncertainty sets $N$ vary. For each uncertainty set $\mathcal{F}_i$ the density parameter $\gamma_i$ is randomly chosen in the interval $(0.09,0.15)$ and the elements of the set $\mathcal{F}_i$ are randomly generated in accordance to the parameter $\gamma_i$. The tables present  the average number of iterations and the computation time in which the algorithm terminates and finds the solution. Table 4 contains the data for the maximisation, while Table 5 presents the results for the minimization problem.

\begin{table}[h!]
\begin{center}
\begin{tabular}{c|c c c}
 $d\setminus N$ & 50  & 100 & 250\\
\hline
25 & 5.5 & 6.2 & 6.3\\
100 & 4.2 & 4.5 & 4.6\\
500 & 4.1 & 4.3 & 4.3\\
2000 & 4.1 & 4.3 & 4.1
\end{tabular}
\caption*{{\footnotesize Table 4a: Average number of iterations for maximization, for non-negative families}}
\end{center}
\end{table}

\begin{table}[h!]
\begin{center}
\begin{tabular}{c|c c c}
 $d\setminus N$ & 50  & 100 & 250\\
\hline
25 & 0.04s & 0.09s & 0.2s\\
100 & 0.12s & 0.25s & 0.72s\\
500 & 0.63s & 1.23s & 3.06s\\
2000 & 3.51s & 6.89s & 187.44s
\end{tabular}
\caption*{{\footnotesize Table 4b: Average computing time for maximization, for non-negative families}}
\end{center}
\end{table}

\begin{table}[h!]
\begin{center}
\begin{tabular}{c|c c c}
 $d\setminus N$ & 50  & 100 & 250\\
\hline
25 & 6.8 & 7.7 & 6.9\\
100 & 5.3 & 4.9 & 5.1\\
500 & 4.2 & 4.1 & 4.6\\
2000 & 4.2 & 4.1 &  4.3
\end{tabular}
\caption*{{\footnotesize Table 5a: Average number of iterations for minimization, for non-negative families}}
\end{center}
\end{table}

\begin{table}[h!]
\begin{center}
\begin{tabular}{c|c c c}
 $d\setminus N$ & 50  & 100 & 250\\
\hline
25 & 0.06s & 0.13s & 0.25s\\
100 & 0.16s & 0.26s & 0.94s\\
500 & 0.75s & 1.18s & 3.01s\\
2000 & 3.31s & 5.77s & 203.28s
\end{tabular}
\caption*{{\footnotesize Table 5b: Average computing time for minimization, for non-negative families}}
\end{center}
\end{table}

Table 6 shows how the number of iterations and computing time vary as the density parameter is changed. The dimension is kept fixed at $d = 600$ and the cardinality of each product set at $\vert\mathcal{F}_i\vert=200$, while we vary the interval $I$ 
(in the first line of the table) from which the density parameter takes value.
For example, the first column of Table 6a shows results for matrices 
that have the ratio of nonzero elements between 0.09 and 0.15. The first line contains the results of the  problem of maximisation o spectral radius (MAX), 
the second does for minimisation  (MIN).  
the problems \\

\begin{table}[h!]
\begin{center}
\begin{tabular}{c|c c c c}
 $I$ & (0.09,0.15)  & (0.16,0.21) & (0.22,0.51) & (0.52,0.76)\\
\hline
MAX & 4.4 & 4.3 & 4.1 & 3.9\\
MIN & 4.5 & 4 & 4.4 & 3.8\\
\end{tabular}
\caption*{{\footnotesize Table 6a: Effects of sparsity, number of iterations
for $d=600, N = 200$}}
\end{center}
\end{table}

\begin{table}[h!]
\begin{center}
\begin{tabular}{c|c c c c}
 $I$ & (0.09,0.15)  & (0.16,0.21) & (0.22,0.51) & (0.52,0.76)\\
\hline
MAX & 2.9s & 2.69s & 2.57s & 2.43s\\
MIN & 2.79s & 2.47s & 3s & 2.35s\\
\end{tabular}
\caption*{{\footnotesize Table 6b: Effects of sparsity, computing time for $d=600, N = 200$}}
\end{center}
\end{table}

\newpage

\begin{center}
\textbf{6.2.  Some conclusions}
\end{center}
\medskip 

The main conclusion from the numerical results is rather surprising. 
The number of iterations seems to depend neither on the dimension, 
nor on the cardinality of the uncertainty sets nor on the sparsity of matrices. 
The usual number of iterations is 3 - 4. The robustness to sparsity can be explained by the fact that  the selective greedy method actually works with 
positive perturbations of matrices. The independence of cardinality 
of the uncertainty sets (in the tables we see that the number of iterations may even decrease in cardinality) is, most likely, explained by the quadratic convergence. 
The distance to the optimal matrix quickly decays to  the tolerance parameter, which is  usually set up between $10^{-8}$ and $10^{-6}$, and the algorithm terminates.  
\medskip

\begin{center}
\textbf{6.3.  Polyhedral product families}
\end{center}
\medskip

Here we test the selective greedy method on product families with  uncertainty sets given by systems of $2d + N$ linear constraints of the form:
\begin{equation*}
\begin{cases}
(\bx,\bb_j)\leqslant 1 & j=1,\ldots,N\\
0\leqslant x_i\leqslant 1 & i=1,\ldots,d
\end{cases}
\end{equation*}
where vectors $\bb_j\in\mathbb{R}^d_+$ are randomly generated and normalized, and $\bx=(x_1,\ldots,x_d)$. As before, $d$ is the dimension and 
for each pair $(d, N)$, ten experiments have been made, the average is put in the table. 
\\

Table 7 shows the test results for  maximizing the spectral radius.

\begin{table}[h!]
\begin{center}
\begin{tabular}{c|c c c}
 $d\setminus N$ & 5  & 10 & 50\\
\hline
10 & 3.6 & 3.6 & 3\\
25 & 4 & 4.2 & 3.8\\
75 & 4.6 & 4.2 & 4.4\\
150 & 4.8 & 4.6 & 4.6 
\end{tabular}
\caption*{{\footnotesize Table 7a: Average number of iterations for polyhedral sets}}
\end{center}
\end{table}

\begin{table}[h!]
\begin{center}
\begin{tabular}{c|c c c}
 $d\setminus N$ & 5  & 10 & 50\\
\hline
10 & 1.42s & 1.46s & 2.35s\\
25 & 4.42s & 5.89s & 15.24s\\
75 & 32.58s & 50.38s & 244.86s\\
150 & 166.43s & 343.33s & 1788.16s 
\end{tabular}
\caption*{{\footnotesize Table 7b: Average computing time for polyhedral sets}}
\end{center}
\end{table}

In polyhedral uncertainty sets the algorithm anyway searches the optimal rows among the 
vertices. However, the number of vertices can be huge. Nevertheless, the 
selective greedy method still needs less than 5 iterations to 
find an optimal matrix. 

\bigskip

\begin{center}
\large{\textbf{7. Applications}}
\end{center}
\bigskip

We consider two possible applications of the greedy method: optimising 
spectral radius of graphs and finding the closest stable matrix. 
\bigskip 

\begin{center}
\textbf{7.2. Optimizing the spectral radius of a graph}
\end{center}
\bigskip

The spectral radius of a graph is the spectral radius of its adjacency matrix $A$.
The problem of maximizing/mimimizing the spectral radius under some restrictions on the graph is well known  in the literature (see~\cite{brualdi, Cve, ESS, friedland, Liu, Olesky} and the references therein).  Some of them deal with product families of 
adjacency matrices and hence can be solved by the greedy method.
For example, to find the maximal spectral radius of a directed graph 
with $d$ vertices and with prescribed numbers of incoming edges  $n_1, \ldots , n_d$.
For this problem, all the uncertainty sets~$\cF_i$ are polyhedral and are given by
the systems of inequalities: 
\begin{equation}\label{eq.graph-max}
 \cF_i \ = \ \Bigl\{x \, \in \, \re^d\ \Bigl| \
\sum_{k=1}^d x_k \, \le \, n_i\, , \, 0\le x_k \le 1, \, k = 1, \ldots , d\,  
\Bigr\}\, .
\end{equation}
The minimal and maximal spectral radii are both attained at extreme points of the uncertainty sets, i.e.,
precisely when the the $i$th row has $n_i$ ones and all other entries are zeros.
If we need to minimize the spectral radius, we define $\cF_i$ in the same way, 
with the only one difference: $\sum_{k=1}^d x_k \, \ge \, n_i$. 

Performing  the greedy method we solve  in each iteration an LP 
(linear programming) problem $(\ba_i , \bv_k) \to \max, \ 
\ba_i \in \cF_i, \ i = 1, \ldots , d$,
 with the sets $\cF_i$ defined in~(\ref{eq.graph-max}). In fact, this LP
 problem is solved explicitly:  the (0,1)-vector $\ba_i$ has 
 ones precisely at positions of the $n_i$ maximal entries of 
 the vector $\bv_k$, and all other components of~$\ba_i$ are zeros. 

Thus, in each iteration, instead of solving $d$ LP problems, which can get computationally demanding even for not so big $d$, we just need to deal with finding the corresponding number of highest entries of the eigenvector $\bv_k$. In Table 8 we present the results of numerical simulations for applying the selective greedy method to this problem, where the number of iteration and the time required for the algorithm to finish are shown. 

\begin{table}[h!]
\begin{center}
\begin{tabular}{c|c c c c}
$d$  & 500 & 1500 & 3000 & 5000\\
\hline
$\#$ & 3 & 3 & 3 & 3\\
$t$ & 0.16s & 0.83s & 3.08s & 7.47s\\
\end{tabular}
\caption*{{\footnotesize Table 8: Row sums are randomly selected integers from the segment [75,100]}}
\end{center}
\end{table}

For example, we apply the greedy method to 
find the maximal spectral radius of a graph with 7 vertices and 
with the numbers of incoming edges  $ (n_1, \ldots , n_7)= (3,2,3,2,4,1,1)$. The algorithm finds  the optimal graph with the adjacency matrix: 

\[
A = 
\left(
\begin{array}{l c c c c c r}
1 {\ } & 0 {\ } & 1 {\ } & 0 {\ } & 1 {\ } & 0 {\ } & 0\\
0 {\ } & 0 {\ } & 1 {\ } & 0 {\ } & 1 {\ } & 0 {\ } & 0\\
1 {\ } & 0 {\ } & 1 {\ } & 0 {\ } & 1 {\ } & 0 {\ } & 0\\
0 {\ } & 0 {\ } & 1 {\ } & 0 {\ } & 1 {\ } & 0 {\ } & 0\\
1 {\ } & 0 {\ } & 1 {\ } & 1 {\ } & 1 {\ } & 0 {\ } & 0\\
0 {\ } & 0 {\ } & 0 {\ } & 0 {\ } & 1 {\ } & 0 {\ } & 0\\
0 {\ } & 0 {\ } & 0 {\ } & 0 {\ } & 1 {\ } & 0 {\ } & 0\\
\end{array}
\right)
\]
with $\rho(A) = 3.21432$.

\bigskip

\begin{center}
\textbf{7.2. Finding the closest stable non-negative matrix}
\end{center}
\bigskip

A matrix is called stable (in the sense of Schur stability) if its 
spectral radius is smaller than one. The {\em matrix stabilization problem}, i.e., the problem of finding the  
closest stable matrix to a given matrix is well known and has been studied in the literature. Finding the closest {\em non-negative} matrix is important 
in applications such as dynamical systems, mathematical economy, population dynamics, etc. (see~\cite{Log, An, NP2, GP2}). Usually the closest matrix is defined in the Euclidean or in the Frobenius norms. In both cases the problem is hard. It was first noted in~\cite{NP2} that in $L_{\infty}$-norm there are efficient methods to 
find the global minimum. Let us recall that 
$\|A\|_{\infty} = \max_{i=1, \ldots , d}\sum_{j=1}^d|a_{ij}|$. 
For a given matrix $A\ge 0$ such that $\rho(A) > 1$, we consider the problem

\begin{equation}\label{eq.stab}
\left\{
\begin{array}{ll}
\|X-A\|_{\infty} \to \min\\
\mbox{subject to} \ X\ge 0, \, \rho(X) \le 1.
\end{array}
\right.
\end{equation}

\medskip

The key observation is that for every positive $r$, the 
ball of radius $r$ in $L_{\infty}$-norm centered in~$A$, i.e.,
the set of non-negative matrices 
$$
\cB(r, A)\ =  \ \Bigl\{X \ge 0\ \Bigl| \ \|X-A\|_{\infty} \le r  \Bigr\}
$$    
is a product family with the polyhedral uncertainty sets 
$\cF_i = \{\bx \in \re^d_+, \ \sum_{j=1}^{d} |x_{j} - a_{ij}| \le r\}$. Hence, for each $r$, the problem 
$\rho(X)\to \min, \ X \in \cB(r, A)$, can be solved by the selective greedy method. 
Then the minimal $r$ such that $\rho(X) \le 1$ is the solution of 
problem~(\ref{eq.stab}). This $r$ can be found merely by bisection.

We found closest stable matrices to randomly generated  matrices $A$ in 
various dimensions up to $d=500$.  The results are given in the Table 9. We remark that for this set of experiments we set the precision of the power method to $10^{-8}$, since for higher dimensions the procedure gets really sensitive due to the fact that the entries of the leading eigenvector get pretty close to one another. This affects the procedure which is dependent on the ordering of the entries of the leading eigenvector.

\begin{table}[h!]
\begin{center}
\begin{tabular}{c|c c c c}
 $d$ & 50  & 100 & 250 & 500\\
\hline
$t$ & 3.4s & 5.38s & 22.98s & 99.44s\\
\end{tabular}
\caption*{{\footnotesize Table 9: Average computing time for finding the closest stable matrix}}
\end{center}
\end{table}

Below is one practical example with an unstable  matrix $A$ of size $d=10$ with $\rho(A) = 9.139125$:

\[
A = 
\left(
\begin{array}{l c c c c c c c c r}
0 {\ } & 0 {\ } & 0  {\ } & 3 {\ } & 5 {\ } & 0 {\ } & 8 {\ }& 0 {\ } & 0 {\ } & 0 \\
8 {\ } & 0 {\ } & 0  {\ } & 0 {\ } & 0 {\ } & 0 {\ } & 0 {\ }& 0 {\ } & 8 {\ } & 0 \\
0 {\ } & 2 {\ } & 0  {\ } & 0 {\ } & 0 {\ } & 4 {\ } & 0 {\ }& 5 {\ } & 0 {\ } & 7 \\
0 {\ } & 0 {\ } & 0  {\ } & 0 {\ } & 0 {\ } & 0 {\ } & 0 {\ }& 0 {\ } & 8 {\ } & 0 \\
1 {\ } & 0 {\ } & 0  {\ } & 0 {\ } & 0 {\ } & 0 {\ } & 0 {\ }& 0 {\ } & 0 {\ } & 0 \\
0 {\ } & 0 {\ } & 0  {\ } & 0 {\ } & 0 {\ } & 7 {\ } & 0 {\ }& 0 {\ } & 0 {\ } & 0 \\
0 {\ } & 0 {\ } & 0  {\ } & 0 {\ } & 6 {\ } & 2 {\ } & 2 {\ }& 0 {\ } & 1 {\ } & 0 \\
0 {\ } & 0 {\ } & 0  {\ } & 0 {\ } & 0 {\ } & 0 {\ } & 1 {\ }& 0 {\ } & 7 {\ } & 0 \\
0 {\ } & 0 {\ } & 0  {\ } & 9 {\ } & 5 {\ } & 0 {\ } & 0 {\ }& 0 {\ } & 1 {\ } & 0 \\
0 {\ } & 0 {\ } & 0  {\ } & 0 {\ } & 0 {\ } & 0 {\ } & 3 {\ }& 4 {\ } & 8 {\ } & 9 \\
\end{array}
\right)
.
\]

The algorithm finds the closest (in $L_{\infty}$-norm) stable matrix:

\[
X = 
\left(
\begin{array}{l l l l l l l l l r}
0 {\ } & 0 {\ } & 0  {\ } & 3 {\ } & 5 {\ } & 0 {\ } & 0.125 {\ }& 0 {\ } & 0 {\ } & 0 \\
0.125 {\ } & 0 {\ } & 0  {\ } & 0 {\ } & 0 {\ } & 0 {\ } & 0 {\ }& 0 {\ } & 8 {\ } & 0 \\
0 {\ } & 2 {\ } & 0  {\ } & 0 {\ } & 0 {\ } & 4 {\ } & 0 {\ }& 5 {\ } & 0 {\ } & 0 \\
0 {\ } & 0 {\ } & 0  {\ } & 0 {\ } & 0 {\ } & 0 {\ } & 0 {\ }& 0 {\ } & 0.125 {\ } & 0 \\
0 {\ } & 0 {\ } & 0  {\ } & 0 {\ } & 0 {\ } & 0 {\ } & 0 {\ }& 0 {\ } & 0 {\ } & 0 \\
0 {\ } & 0 {\ } & 0  {\ } & 0 {\ } & 0 {\ } & 0 {\ } & 0 {\ }& 0 {\ } & 0 {\ } & 0 \\
0 {\ } & 0 {\ } & 0  {\ } & 0 {\ } & 6 {\ } & 2 {\ } & 2 {\ }& 0 {\ } & 1 {\ } & 0 \\
0 {\ } & 0 {\ } & 0  {\ } & 0 {\ } & 0 {\ } & 0 {\ } & 1 {\ }& 0 {\ } & 7 {\ } & 0 \\
0 {\ } & 0 {\ } & 0  {\ } & 2.125 {\ } & 5 {\ } & 0 {\ } & 0 {\ }& 0 {\ } & 0 {\ } & 0 \\
0 {\ } & 0 {\ } & 0  {\ } & 0 {\ } & 0 {\ } & 0 {\ } & 3 {\ }& 4 {\ } & 8 {\ } & 1 \\
\end{array}
\right)
.
\]

\bigskip

\bigskip

\bigskip

\begin{center}
\large{\textbf{8. Implementation details of the selective greedy method}}
\end{center}
\bigskip

\begin{center}
\textbf{8.1. Rounding errors. The tolerance parameter.}
\end{center}
\bigskip

Even though the selective greedy method, in theory, can safely and effectively be used on non-negative product families, the cycling might occur due to the computational errors, especially in the case of very sparse  matrices. The following example sheds some light on this issue. 

We run the selective greedy algorithm to the product family $\mathcal{F} = \mathcal{F}_1\times\cdots\times\mathcal{F}_7$ of $(0,1)$-matrices of dimension $7$, for solving the minimization problem. Each uncertainty set $\mathcal{F}_i$ consists of $(0,1)$-vectors containing from $1$ to $4$ ones.   The selective greedy algorithm cycles between the following two matrices:
\[ A \ = \  
\begin{pmatrix}
0 & 1 & 1 & 0 & 0 & 1\\
0 & 1 & 0 & 0 & 0 & 1\\
0 & 1 & 0 & 0 & 1 & 0\\
0 & 1 & 1 & 0 & 1 & 1\\
0 & 1 & 0 & 0 & 1 & 0\\
0 & 1 & 0 & 0 & 0 & 1
\end{pmatrix}
\quad
\rightleftarrows
\quad
 A' \ = \  
\begin{pmatrix}
0 & 1 & 0 & 0 & 1 & 1\\
0 & 1 & 0 & 0 & 0 & 1\\
0 & 0 & 0 & 0 & 1 & 1\\
0 & 1 & 1 & 0 & 1 & 1\\
0 & 0 & 1 & 0 & 0 & 1\\
0 & 0 & 1 & 0 & 0 & 1
\end{pmatrix}
\]

Both matrices have the same leading eigenvalue $\lambda_{max} = 1$, and the algorithm computes the same leading eigenvector 
$$\bv = (0.00106389,\ 0.14841278,\\ 0.73568055,\ 0.44311056,\ 0.14841278,\ 0.14841278,\ 0.44311056).$$ 
Both matrices are minimal in the second, fifth, and sixth rows (the rows that get changed) with the respect to the vector $\bv$. In addition, we have $(\ba_i,\bv) = (\ba'_i,\bv)$ for  $i\in\{2,5,6\}$. According to the algorithm, those rows should remain unchanged, and algorithm have actually to  terminate  with the matrix~$A$. The  explanation for this contradictory behaviour is that the machine, due to the calculation imprecisions, keeps miscalculating the above dot products, first taking the ``problematic'' rows $\ba'_i$ of the second matrix as the optimal, and then rolling back to the rows $\ba_i$, seeing them as optimal.\\

One of the most straightforward strategy to resolve this issue is to simply round all calculation to some given number of decimals. This approach is particularly efficient if we are not concerned with high precision. However, if we want our calculations to have a higher order of precision, we need to resort to other strategies.\\

Another idea to deal with this issue is to modify the part of the algorithm that dictates the conditions under which the row will change. 
Here we will use notation for the maximization, although for everything is completely analogous for the minimization case.
 
Let $\ba_i^{(k)}$ be the $i$-th row of a matrix $A_k$ obtained in the $k$th iteration, $\bv_k$ its leading eigenvector, and $\ba_i^{(k+1)} = \max_{a\in\mathcal{F}_i}(\ba,\bv_k)$, different from $\ba_i^{(k)}$. 
Unless the dot computed products $(\ba_i^{(k)},v_k)$ and $(\ba_i^{(k+1)},v_k)$ are the same, the row $\ba_i^{(k)}$ will get replaced by  $\ba_i^{(k+1)}$.
However, as can be seen from the examples above, the change may occur even if those dot products are the same. 
The undesired change may also occur even if the vector $\ba_i^{(k+1)}$ is not truly optimal, but computed dot products are really close to each other, which will lead to rolling back to the vector $\ba_i^{(k)}$ in the next iteration, and thus cycling.
To counter this, we impose the rule that the row $\ba_i^{(k)}$ will get replaced by the row $\ba_i^{(k+1)}$ if and only if $\ba_i^{(k+1)} = \max_{a\in\mathcal{F}_i}(\ba,\bv_k)$ and $\vert(\ba_i^{(k)},\bv_k) - (\ba_i^{(k+1)},\bv_k)\vert\geqslant\delta$, where $\delta$ is a small tolerance parameter.
In other words, in the $(k+1)$st iteration, the row $\ba_i^{(k)}$ will not  be replaced by the new row $\ba_i^{(k+1)}$ if their computed scalar products with the vector $\bv_k$ are really close to each other. This modification of the algorithm converges to a point of approximate maximum/minimum. The proof of convergence 
along with the estimates for the rates of convergence from Section 5 can be adopted to this modification in a straightforward manner.

\bigskip

\begin{center}
\textbf{8.2. Dealing with reducible families}
\end{center}
\bigskip

As stated in subsection 2.2, positivity of the leading eigenvector $\bv_k$ obtained in the final iteration of the selective greedy algorithm for maximization problem is a guarantee that the obtained solution is correct. However, if  $\bv_k$ contains some zero entries, our computed solution might not be the optimal one. This is not an issue when working with the irreducible product families, since in that case the leading eigenvector of computed matrix is always positive~\cite[Lemma~4]{P16}. 

One way to resolve this issue is to compute the Frobenius factorization 
of the family~$\cF$, after which the  maximization problem will be reduced to 
several similar problems in smaller dimensions (see~~\cite[Section~3.3]{P16}). 
The algorithm of Frobenius factorization  is fast, it can be found in~\cite{T}. 
In practice it is often possible to manage without the Frobenius factorization,   
since having only one irreducible matrix is enough to render the whole family irreducible. A simple strategy when working with reducible families is to simply include a cyclic permutation matrix $P$, which is irreducible, multiplied by a small parameter $\alpha$. 

\bigskip

\begin{center}
\textbf{Acknowledgements}. 
\end{center}
\bigskip 

The authors are grateful to 
two anonymous Referees 
for attentive reading and for many useful suggestions.


\bigskip 

\begin{thebibliography}{NN}



\bibitem{Akian1}
M.Akian, S.Gaubert, J.Grand-Clément, and J.Guillaud, 
\newblock {\em  The operator approach to entropy games}, 
\newblock  Theory Comput. Syst. (2019), published electronically: https://doi.org/10.1007/s00224-019-09925-z
\smallskip 


\bibitem{Al}
E.Altman, 
\newblock {\em Constrained Markov decision process}, 
\newblock INRIA, 2004. 
\smallskip

\bibitem{An}
J.Anderson,
\newblock {\em Distance to the nearest stable Metzler matrix},
\newblock (2017), arXiv:1709.02461v1 
\smallskip


\bibitem{Asarin}
E.Asarin, J.Cervelle, A.Degorre, C.Dima, F.Horn, and V.Kozyakin, 
\newblock {\em   Entropy games
and matrix multiplication games}, 
\newblock  In 33rd Symposium on Theoretical Aspects of Computer
Science (STACS 2016), February 2016, Orl\'eans, France, 11:1 -– 11:14. 
\smallskip 

\bibitem{A}
K.E.Atkinson, 
\newblock {\em An introduction to numerical analysis},
\newblock  John Wiley \& Sons (1989). 
\smallskip 

\bibitem{BN}
V.D.Blondel and Y.Nesterov,
\newblock {\em Polynomial-time computation of
the joint spectral radius for some sets of nonnegative
matrices},
\newblock  SIAM J. Matrix Anal. Appl. 31 (2009), no 3,
865--876.
\smallskip

\bibitem{BP}
A.Berman and R.J.Plemmons,
\newblock {\em Nonnegative matrices in the mathematical sciences},
\newblock  Acedemic Press, Now York, 1979.
\smallskip

\bibitem{brualdi}
R.A.Brualdi and  A.J.Hoffman,
\newblock {\em  On the spectral radius of
$(0,1)$-matrices},
\newblock Linear Alg. Appl., 65 (1985),
133--146.
\smallskip

\bibitem{CL}
D.G.Cantor and S.A.Lippman,
\newblock {\em Optimal investment selection with a multiple of projects},
\newblock Econometrica, 63 (1995),  1231--1240.
\smallskip

\bibitem{CJB}
B.C.Csáji, R.M.Jungers, and V.D.Blondel,  PageRank optimization by edge selection, 
\newblock  {\em Discrete Applied Mathematics}, 169 (2014), 73--87.
\smallskip 

\bibitem{Cve}
D.Cvetkovi\'c, M.Doob, and H.Sachs,
\newblock {\em Spectra of
graphs. Theory and application},
\newblock Academic Press, New York (1980).
\smallskip

\bibitem{ESS}
G.M.Engel, H.Schneider, and S.Sergeev,
\newblock {\em On sets of eigenvalues of matrices with prescribed row sums and prescribed graph},
\newblock  Linear Alg. Appl., 455 (2014), 187-209
\smallskip


\bibitem{FM}
L.Fainshil and M.Margaliot,
 \newblock {\em A maximum principle for the stability analysis of positive bilinear control systems with applications to positive linear switched systems},
 \newblock SIAM J. Control Optim.  50  (2012),  no. 4, 2193--2215.
\smallskip

\bibitem{FABG}
O.Fercoq,  M.Akian,  M.Bouhtou, and  S.Gaubert,  
\newblock {\em Ergodic control and polyhedral approaches to PageRank optimization},  \newblock IEEE Transactions on Automatic Control, 58 (2013), no 1, 134--148.
\smallskip 

\bibitem{friedland}
S.Friedland,
\newblock {\em The maximal eigenvalue of $0$-$1$ matrices
with prescribed number of ones},
\newblock  Linear Alg. Appl.,
69 (1985), 33--69.
\smallskip




\bibitem{G}
F.R.Gantmacher, 
\newblock {\em The theory of matrices},  
\newblock Chelsea, New York, 2013.
\smallskip 

\bibitem{GG}
V.Goyal, J.Grand-Clement, 
\newblock {\em A first-order approach to accelerated value iteration},  
\newblock 	arXiv:1905.09963
\smallskip 

\bibitem{GP2}
N.Guglielmi and V.Protasov, 
\newblock {\em On the closest stable/unstable nonnegative matrix
and related stability radii},  
\newblock SIAM J. Matrix Anal. 39 (2018), no 4,  1642--1669. 
\smallskip 

\bibitem{JPB}
R.Jungers, V.Yu.Protasov,  and V.Blondel,
\newblock  {\it Efficient algorithms for deciding the type of
 growth  of products of integer matrices},
\newblock Linear Alg.  Appl.,  428 (2008),  no 10,  2296--2312.
\smallskip

\bibitem{HMZ} 
T.D. Hansen, P.B. Miltersen, and U. Zwick, 
\newblock {\em Strategy iteration is strongly polynomial for 2-player turn-based stochastic
games with a constant discount factor}, 
\newblock  In Innovations in Computer Science 2011, pages 253–-263. 
Tsinghua University Press,
2011
\smallskip 

\bibitem{HK}
A.J.Hoffman and R.M.Karp, 
\newblock {\em  On nonterminating stochastic games}, 
\newblock  Management Science.
Journal of the Institute of Management Science. Application and Theory Series, 12 (1966), 359-–370.
\smallskip 

\bibitem{HJ}
R.A.Horn and C.R.Johnson,
\newblock {\em Matrix analysis},
\newblock  Cambridge University Press, 1990.
\smallskip

\bibitem{Koz1}
V.S.Kozyakin,
\newblock {\em A short introduction to asynchronous systems},
\newblock in Proceedings of the Sixth International
Conference on Difference Equations (Augsburg, Germany 2001): New Progress
in Difference Equations, B. Aulbach, S. Elaydi, and G. Ladas, eds., CRC Press, Boca
Raton, FL (2004),  153--166.
\smallskip

\bibitem{Koz2}
V.S.Kozyakin,
\newblock {\em 
Hourglass alternative and the finiteness conjecture for the spectral characteristics of sets of non-negative matrices}, 
\newblock Lin. Alg. Appl.  489 (2016), 167--185. 
\smallskip  


\bibitem{LA}
H.\,Lin and P.J.\,Antsaklis,
\newblock {\em Stability and stabilizability of switched linear
systems: a survey of recent results},
\newblock IEEE Trans. Autom. Contr., 54 (2009), no 2, 308--322.
\smallskip

\bibitem{Liu}
B.Liu,
\newblock {\em  On an upper bound of the spectral radius of graphs},
\newblock Discrete Math., 308 (2008), no 2, 5317--5324.
\smallskip

\bibitem{Log}
D.O.\,Logofet,
\newblock {\em Matrices and graphs: stability problems in mathematical ecology},
\newblock CRC Press, Boca Raton, 1993.
\smallskip

\bibitem{M}
K.S.Miller,
\newblock  {\em Linear difference equations},
\newblock Elsevier, Amsterdam, 2000.
\smallskip

\bibitem{MS}
 Y.\,Mansour and S.\,Singh, 
\newblock {\em On the complexity of policy iteration},  
\newblock In Proceedings of the 15th International Conference on Uncertainty in AI, 
(1999),  401--408.
\smallskip


\bibitem{NN}
Y.\,Nesterov and 
A.\,Nemirovski, 
\newblock  {\em Finding the stationary states of Markov chains by iterative
methods},  
 \newblock Appl. Math. Comput. 
 255 (2015),  58--65. 
\smallskip 

\bibitem{NP1}
Y.\,Nesterov, V.Yu.\,Protasov,
\newblock {\em Optimizing the spectral radius},
\newblock SIAM J. Matrix Anal. Appl. 34 (2013), no 3,
999--1013.
\smallskip

\bibitem{NP2}
Y.Nesterov, V.Yu.Protasov,
\newblock {\em Computing closest stable
non-negative matrix}, 
\newblock SIAM J. Matrix. Anal. Appl. no. 1, 1-–28, (2020), 

\smallskip

\bibitem{Olesky}
D.D.Olesky, A.Roy, and P. van den Driessche,
\newblock {\em  Maximal
graphs and graphs with maximal spectral radius},
\newblock Linear Alg. Appl., 346 (2002), 109--130.
\smallskip

\bibitem{Par}
B.G.Pachpatte,
\newblock  {\em Integral and finite difference
inequalities and applications},
\newblock W. A. Benjamin, Inc., New
York-Amsterdam, 1968.
\smallskip

\bibitem{PY}
I.Post, Y.Ye, 
\newblock {\em The simplex method is strongly polynomial for deterministic markov decision processes}, 
\newblock Mathematics of Operations Research,  published online (2015), 
https://doi.org/10.1287/moor.2014.0699

\smallskip 

\bibitem{P}
M.L.\,Putermam, 
\newblock Markov Decision Processes, 
\newblock  Wiley (1994). 
\smallskip 


\bibitem{P16}
V.Yu.Protasov,
\newblock {\em The spectral simplex method},
\newblock Math. Prog., 156 (2016), 485--511
\smallskip

\bibitem{P19}
V.Yu.Protasov,
\newblock {\em How to make the Perron eigenvector simple},
\newblock CALCOLO, 56 (2019) no 17. 
\smallskip

\bibitem{R}
U.G.Rothblum, 
\newblock {\em Multiplicative markov decision chains}, 
 \newblock Math. Oper. Research,
9(1984), 6-–24.

\bibitem{S}
V.Schvatal
\newblock {\em Linear programming},
\newblock W. H. Freeman, 1983.
\smallskip


\bibitem{SS}
G.W.Stewart and J.G.Sun,
\newblock {\em Matrix perturbation theory},
\newblock  Academic Press, New York, 1990.
\smallskip

\bibitem{T}
R.Tarjan, 
\newblock {\em Depth first search and linear graph algorithms}, 
\newblock SIAM J. Comput. 1 (2) (1972),  146-–160.
\smallskip 

\bibitem{as}
A. G. Vladimirov, N. Grechishkina, V. Kozyakin, N. Kuznetsov, A. Pokrovskii, and D.
Rachinskii,
\newblock {\em Asynchronous systems: Theory and practice},
\newblock Inform. Processes, 11 (2011), 1--45.
\smallskip 

\bibitem{Y}
 Y.\,Ye. 
\newblock {\em The simplex and policy-iteration methods are strongly polynomial for the markov decision problem with a fixed
discount rate}, 
\newblock  Mathematics of Operations Research, 36 (2011), no 4, 593–-603.
\smallskip 

\end{thebibliography}
\end{document}